\documentclass[12pt]{article}
\usepackage{amsmath}
\usepackage{amssymb}
\usepackage{amsthm}

\usepackage[margin=2.0cm]{geometry}
\usepackage{hyperref}
\usepackage{multirow}

\usepackage{color}
\usepackage{enumerate}
\usepackage{tikz,lmodern}
\usepackage{amsfonts}
\usetikzlibrary{arrows, decorations.markings}

\newtheorem{lem}{Lemma}[section]%
\newtheorem{theorem}[lem]{Theorem}%
\newtheorem{defi}[lem]{Definition}%
\newtheorem{cor}[lem]{Corollary}%
\newtheorem{prop}[lem]{Proposition}%
\def\a{\alpha}
\def\b{\beta}
\def\g{\gamma}\def\d{\delta}

\def\s{\sigma}

\def\G{\Gamma}
\def\D{\bf D}

\def\nd{\mathrel{\bigm|\kern-.7em/}}

\def\f{\noindent}

\def\Aut{\hbox{\rm Aut}}

\def\Ker{\hbox{\rm Ker}}

\def\Cay{\hbox{\rm Cay}}

\def\Inn{\hbox{\rm Inn}}
\def\Arc{\hbox{\rm Arc}}

\def\demo{\f {\bf Proof.}\hskip10pt}

\def\mz{{\mathbb Z}}

\def\D{\hbox{\rm D}}

\def\S{{\rm S}}

\def\H{\mathcal{H}}
\def\G{\mathcal {G}}

\numberwithin{equation}{section}

\tikzstyle{vecArrow}=[semithick, decoration={markings,mark=at position 1 with {\arrow[semithick]{open triangle 60}}},
   double distance=1.4pt, shorten >= 5.5pt,
   preaction = {decorate},
   postaction = {draw,line width=1.4pt, white,shorten >= 4.5pt}]
\begin{document}
\title{On isomorphisms of $m$-Cayley digraphs}
\author{\\Xing Zhang, Yuan-Quan Feng\footnotemark, Fu-Gang Yin, Jin-Xin Zhou\\
{\small\em School of mathematics and statistics, Beijing Jiaotong University, Beijing, 100044, P.R. China}\\
}

\renewcommand{\thefootnote}{\fnsymbol{footnote}}
\footnotetext[1]{Corresponding author.
E-mails:
21118006@bjtu.edu.cn (X. Zhang),
yqfeng@bjtu.edu.cn (Y.-Q. Feng),
18118010@bjtu.edu.cn (F.-G. Yin),
jxzhou@bjtu.edu.cn (J.-X. Zhou)}

\date{}
\maketitle

\begin{abstract}
The isomorphism problem for digraphs is a fundamental problem in graph theory. This problem for Cayley digraphs has been extensively investigated over the last half a century. In this paper, we consider this problem for $m$-Cayley digraphs which are generalization of Cayley digraphs. Let $m$ be a positive integer. A digraph admitting a group $G$ of automorphisms acting semiregularly on the vertices with exactly $m$ orbits is called an {\em $m$-Cayley digraph} of $G$. In particular, $1$-Cayley digraph is just the Cayley digraph. We first characterize the normalizer of $G$ in the full automorphism group of an $m$-Cayley digraph of a finite group $G$. This generalizes a similar result for Cayley digraph achieved by Godsil in 1981. Then we use this to study the isomorphisms of $m$-Cayley digraphs. The CI-property of a Cayley digraph (CI stands for `Cayley isomorphism') and the DCI-groups (whose Cayley digraphs are all CI-digraphs) are two key topics in the study of isomorphisms of Cayley digraphs. We generalize these concepts into $m$-Cayley digraphs by defining $m$CI- and $m$PCI-digraphs, and correspondingly, $m$DCI- and $m$PDCI-groups. Analogues to Babai's criterion for CI-digraphs are given for $m$CI- and $m$PCI-digraphs, respectively. With these we then classify finite $m$DCI-groups for each $m\geq 2$, and finite $m$PDCI-groups for each $m\geq 4$. Similar results are also obtained for $m$-Cayley graphs. Note that 1DCI-groups are just DCI-groups, and the classification of finite DCI-groups is a long-standing open problem that has been worked on a lot.

\bigskip
\noindent {\bf Keywords:}  $m$-Cayley digraph, $m$-PCayley digraph, Cayley isomorphism, semiregular group. \\
{\bf 2010 Mathematics Subject Classification:} 05C25, 20B25.
\end{abstract}

\section{Introduction\label{Sec1}}

A digraph $\Gamma$ is  an ordered pair $(V(\Gamma),  \mathrm{Arc}(\Gamma))$, where $V(\Gamma)$ is a non-empty set and  $\mathrm{Arc}(\Gamma)$, called the arc set of $\Gamma$, is a subset of $V(\Gamma)\times V(\Gamma)$. The digraph $\Gamma$ is called a graph if  $\mathrm{Arc}(\Gamma)$ is symmetric, that is, $(u,v)\in \mathrm{Arc}(\Gamma)$ implies $(v,u)\in \mathrm{Arc}(\Gamma)$, where $\Gamma$ has the edge set $E(\Gamma)=\{\{u,v\}\ :\ (u,v)\in \mathrm{Arc}(\Gamma)\}$. On the other way, a graph $\Sigma$ is identified as a digraph by defining $\mathrm{Arc}(\Sigma)=\{(u,v)\ :\ \{u,v\}\in E(\Sigma)\}$. {\em The complement $\Gamma^{c}$} of a digraph $\Gamma$ is the digraph with vertex set $V(\Gamma)$ and arc set $\Arc(\Gamma^{c})=\{(u,v)\ :\ (u,v)\not\in\Arc(\Gamma)\}$. Denote by $\Aut(\Gamma)$ the automorphism group of $\Gamma$. Then $(\Gamma^c)^c=\Gamma$ and $\Aut(\Gamma)=\Aut(\Gamma^c)$. Clearly, one of $\Gamma$ and $\Gamma^c$ must be connected. Let $U$ and $V$ be subsets of $V(\Gamma)$. Denote by $[U]_\Gamma$ the induced subdigraph of $\Gamma$ with vertex set $U$ and arc set $\mathrm{Arc}([U]_\Gamma)=\{(u_1,u_2)\in \mathrm{Arc}(\Gamma)\ :\ u_1,u_2\in U\}$.
Denote by $[U,V]_\Gamma$ the subdigraph of $\Gamma$ with vertex set $U\cup V$ and arc set $\mathrm{Arc}([U,V]_\Gamma)=\{(u,v)\in \mathrm{Arc}(\Gamma)\ :\ u\in U, v\in V\}$. If there is no confusion, we omit the subscript $\Gamma$ for $[U]_\Gamma$ and $[U,V]_\Gamma$. Throughout this paper, all digraphs are finite and simple (no self-loop or  multi-arc), and all groups are finite.

Let $\Gamma$ be a digraph  and let $G\leq \Aut(\Gamma)$ be semiregular on $V(\Gamma)$ with $m$ orbits (semiregular means that the stabilizer $G_{u}=1$ for every $u\in V(\Gamma)$). Then $\Gamma$ can be presented as {\em an $m$-Cayley digraph} $\Gamma=\Cay(G,S_{i,j}: 1\leq i,j\leq m)$ of $G$ (see \cite{DFS}) with respect to $m^2$ subsets $S_{ij}$ of $G$  such that
\begin{align*}
V(\Gamma)&=\bigcup_{1\leq i\leq m}  G_i, \text{ where } G_i=\{x_i :  x\in G\},\\
\mathrm{Arc} (\Gamma)&=\bigcup_{1\leq i,j\leq m} \{(x_i, (sx)_j) :  s\in S_{i,j},x\in G\}.
\end{align*}
In particular, we may identify the semiregular automorphism group $G$ with the right multiplication of $G$, namely $R_m(G)=\{R_m(g):g \in G \}$, where $R_m(g)$ is the permutation on $V(\Gamma)$ defined by:
\[
R_m(g): x_i \mapsto (xg)_i, \text{ for all } x  \in G
\text{ and } 1\leq i \leq m.\]
If $S_{ii}=\emptyset$ for all $i$, then $\Gamma$ is called an \emph{$m$-PCayley digraph}.
Note that $\Gamma$ is an $m$-Cayley graph if and only if $S_{ij}=S_{ji}^{-1}$ for each pair $(i,j)$, and the identity $1$ of $G$ is not in each $S_{ii}$ because $\Gamma$ has no self-loop.
The $m$-Cayley digraphs contain several well-known families of digraphs, such as Cayley digraphs ($1$-Cayley digraphs), bi-Cayley graphs ($2$-Cayley  digraphs), and Haar digraphs  ($2$-PCayley digraphs) in the literature.

A Cayley (di)graph $\mathrm{Cay}(G, S)$ is said to be \emph{CI} (CI stands for Cayley isomorphism) if, for any
Cayley (di)graph $\mathrm{Cay}(G, T)$, whenever $\mathrm{Cay}(G, S)\cong \mathrm{Cay}(G, T)$ we have $S^\alpha=T$ for some $\alpha \in \mathrm{Aut}(G)$. A group $G$ is called  a \emph{DCI-group} or \emph{CI-group} if any Cayley digraph or graph of $G$ is CI, respectively. There are two long open questions for CI-problems:
\begin{itemize}
\item Which Cayley (di)graphs for a group $G$ are CI?
\item Which groups are DCI-groups or CI-groups?
\end{itemize}

Numerous papers have been published on these two questions, for which the reader may refer to the survey papers \cite{Li,XuM}. One of the remarkable achievements on these questions is the classifications of cyclic DCI-groups and cyclic CI-groups, the study of which dates back to 1967 when \'Ad\'am conjectured in \cite{Ad} that every cyclic group is a CI-group. In 1970, Elspas and Turner~\cite{Elspas} disproved \'Ad\'am conjecture. Since then, a lot of work have been devoted to the classification of cyclic CI-groups. Based on contributions of many researchers like Elspas and Turner~\cite{Elspas}, Djokovi\'c~\cite{Djokovic}, Turner~\cite{Turner}, Babai~\cite{Babai}, Alspach and Parsons~\cite{Alspach}, Godsil~\cite{Godsil} and P\'alfy~\cite{Palfy}, the cyclic DCI-groups and CI-groups were finally classified by Muzychuk~\cite{Muzychuk,M.Muzychuk} in 1997. After this, many people have worked on seeking general DCI-groups and CI-groups, and a lot of work have been done over last fifty years (see for instance \cite{AN,BDM,D2,DE1,FK,Kov,Mu3,SM,N,So,Sp,Sp1,Xie1,C.H.Li9}). In 2007, Li et al. in \cite{C.H.Li8} gave an explicit list of candidates for finite CI-groups. However, it is very hard to determine whether a particular group in the list is a DCI-group or a CI-group. For example, the dihedral group is a group in the list, and the problem of classifying dihedral CI-groups is still widely open (see \cite{Babai,conder,DE1,DT,Xie2,Xie3,Kovacs} for some progress on this problem).

In 2015, Arezoomand and Taeri~\cite{Arez} initiated the study of isomorphism problem for $2$-Cayley graphs, and they introduced the following concept: a $2$-Cayley graph $\Gamma=\mathrm{Cay}(G,\{S_{11},S_{12},S_{21},S_{22}\})$ is called a \emph{SCI-graph} (SCI stands for semi-Cayley isomorphism) if, for any $2$-Cayley graph $\Sigma$ of $G$ isomorphic to $\Gamma$, there exists some $\a\in N_{\S_{V(\Gamma)}}(R(G))$ such that $\Gamma^\a=\Sigma$ (see \cite[Lemma~1.1]{Arez} and Proposition~\ref{c-R(G)} for the definition of such an $\a$). In particular,  $\Gamma$ is called a \emph{$0$-type-SCI-graph} if $\Gamma$ is a $2$-PCayley graph that satisfies the above conditions.

In this paper, we shall study the isomorphism problem for the general $m$-Cayley digraphs. Similar to the concepts of CI-(di)graphs and SCI-graphs, we introduce the $m$CI-(di)graphs.

\begin{itemize}
\item An $m$-Cayley (di)graph $\Gamma$ of $G$ is said to be \emph{mCI} (mCI stands for $m$-Cayley isomorphism) if, for any $m$-Cayley (di)graph $\Sigma$ of $G$ isomorphic to $\Gamma$, there exists some $\a\in  N_{\mathrm{S}_{V(\Gamma)}}(R(G))$ such that $\Gamma^\a=\Sigma$.
\end{itemize}

Similarly, for $m$-PCayley (di)graphs, we define $m$PCI-(di)graphs.

\begin{itemize}
\item  An $m$-PCayley (di)graph $\Gamma$ of $G$ is said to be \emph{mPCI} if, for any $m$-PCayley (di)graph $\Sigma $ of $G$  such that there is an isomorphism from $\Gamma$ to $\Sigma$ keeping $\{G_1,G_2,\cdots,G_m\}$ invariant, there is some $n \in  N_{\mathrm{S}_{V(\Gamma)}}(R(G))$ such that $\Gamma^n=\Sigma$.
\end{itemize}

Note that 1CI-(di)graph is exactly the same as CI-(di)graph, 2PCI-graph is the same as $0$-type-SCI-graph in~\cite{Arez}, and by Corollary~\ref{c-R(G)} and Proposition~\ref{mCImPCI}, we shall see that 2CI-graph is exactly the same as SCI-graph.

In the study of isomorphism problem for (di)graphs, it is quite helpful to gain as much information as possible for the group of automorphisms of the (di)graph. This is particular true for Cayley (di)graphs. Let $\Gamma$ be a Cayley (di)graph of a group $G$. In 1981, Godsil~\cite{Godsil1} determined the normalizer of $R(G)$ in $\Aut(\Gamma)$. This result has been successfully used in the study of automorphisms and isomorphisms of Cayley digraphs (see for example \cite{Godsil1,Godsil1983,Li,XuM}). Another important tool in the study of isomorphisms of Cayley (di)graphs is Babai's criterion for CI-digraph: A Cayley (di)graph $\Gamma=\mathrm{Cay}(G,S)$ of a group $G$ is CI if and only if all regular subgroups of $\mathrm{Aut}(\Gamma)$ isomorphic to $G$ are conjugate. This was proved by Babai in \cite{Babai-1977}, and plays a key role in the study of CI-digraphs. Considering the importance of the above two results, it is natural to generalize them to the general $m$-Cayley (di)graphs.

Let $m$ be a positive integer and let $\Gamma$ be an $m$-Cayley (di)graph of a group $G$. When $m=1$, it is well-known that the normalizer $N_{\mathrm{S}_{V(\Gamma)}}(R(G))$ of $R(G)$ in $\mathrm{S}_{V(\Gamma)}$ equals $R(G)\rtimes\Aut(G)$, which is called the {\em holomorph} of $G$ (see~\cite[pp.36-37]{DRobinson}). In this paper, we first generalize this result to all $R_m(G)$ by determining the normalizer of $R_m(G)$ in the symmetric group $\mathrm{S}_{V(\Gamma)}$ for all $m$ (see Theorem~\ref{thN}). Applying this, we determine the normalizer of $R_m(G)$ in $\Aut(\Gamma)$ and some subgroups of the normalizer (see Theorem~\ref{thNmC}), which generalizes the above mentioned result of Godsil for Cayley (di)graphs in \cite{Godsil1}. Then we generalize the Babai's criterion for CI-(di)graphs to all $m$CI-(di)graphs and $m$PCI-(di)graphs by using the characterization of the normalizer of $R_m(G)$ in $\Aut(\Gamma)$ (see Theorems~\ref{BabaiSimilar} and \ref{PCI-BabaiSimilar}).
The criterions for $m$CI-(di)graphs and $m$PCI-(di)graphs are then used to classify the $m$DCI-groups, $m$CI-groups, $m$PDCI-groups and $m$PCI-groups which are defined below.

A group $G$ is called an \emph{$m$DCI}-group or an \emph{$m$CI}-group if every $m$-Cayley digraph or graph of $G$ is $m$CI, respectively.
Similarly, we can define $m$PDCI-groups and $m$PCI-groups. For $m\geq 2$, we shall prove that every $m$DCI-group ($m$CI-group, $m$PDCI-group and $m$PCI-groups, respectively) is $(m-1)$DCI-group ($(m-1)$CI-group, $(m-1)$PDCI-group and $(m-1)$PCI-groups, respectively) (see Theorem~\ref{mtom-1}). This together with \cite[Theorem B]{Arez} enable us to obtain a classification of $m$DCI-groups and $m$CI-groups for every $m\geq 2$ (see Theorem~\ref{mCI-groups}):
\begin{itemize}
\item A group $G$ is $m$CI if and only if either $m=2$ and $G=1$ or $\mz_3$, or $m\geq 3$ and $G=1$;
\item A group $G$ is $m$DCI if and only if $G=1$.
\end{itemize}
Furthermore, we also give a classification of $m$PCI-groups and $m$PDCI-groups for every $m\geq 4$ (see Theorem~\ref{mPCI-groups}):
\begin{itemize}
\item A group $G$ is $m$PCI if and only if $G=1$, $\mathbb{Z}_2$, or $G=\mathbb{Z}_3$ or $D_6$ with $m=4$ or $5$.
\item A group $G$ is $m$PDCI if and only if $G=1$ or $\mathbb{Z}_2$.
\end{itemize}
It is worthy to note that a subgroup of an $m$PCI-group with $m\geq 3$ is also $m$PCI, see Theorem~\ref{subgroupsmPCI}. However, we do not know whether it is true for $m=2$. It is well known that a subgroup of a CI-group is CI (see \cite[Lemma 3.2]{Babai}).

\medskip
The paper is organised as follows. In Section~\ref{Sec2}, we determine the normalizer $N_{\S_{\Omega}}(G)$ of a semiregular permutation group $G$ on a set $\Omega$ in the symmetric $\S_{\Omega}$, the centralizer $C_{\S_{\Omega}}(G)$ of $G$ in $\S_{\Omega}$, and the kernel of $N_{\S_{\Omega}}(G)$ acting on the orbit set of $G$. In Section~\ref{sec3}, based on the results in Section~\ref{Sec2}, we characterize the normalizer of $R(G)$ in the automorphism group of an $m$-Cayley digraph of a group $G$ and some subgroups of the normalizer. In Section~\ref{sec4}, we prove the criterions for $m$CI-(di)graphs and $m$PCI-(di)graphs, and in Section~\ref{sec5}, we investigate the classifications of $m$DCI-groups, $m$CI-groups, $m$PCI-groups and $m$PDCI-groups.

To end this section we fix some notation used in this paper. The notation for groups in this paper is standard; see \cite{Atlas} for example. In particular, denote by $\mz_{m}$ and $\mz_m^*$ the additive group of integer numbers modulo $m$, and the multiplicative group of numbers coprime to $m$ in $\mz_m$, by $\D_{m}$ the dihedral group of order $m$, and by $\S_m$ the symmetric group of degree $m$, respectively. Also, we denote by $\S_{\Omega}$ the symmetric group on the set $\Omega$.
For a prime $p$,  we use $\mz_p^m$ to denote the elementary abelian group of order $p^m$. For two groups $A$ and $B$, $A\times B$ stands for the direct product of $A$ and $B$, and $A\rtimes B$ for a split extension or a semi-direct product of $A$ by $B$.

\section{Normalizer of a semiregular group\label{Sec2}}

In this section we give an explicit construction of the normalizer of a semiregular permutation group on a set $\Omega$ in the symmetric group $\S_\Omega$, and study some subgroups of the normalizer. Recall that two permutation groups $X \leq \mathrm{S}_{\Delta}$ and $X'\leq \mathrm{S}_{\Delta'}$ are called \emph{permutation isomorphic} if there exists a bijection $\lambda: \Delta \mapsto \Delta'$ and a group isomorphism $\varphi: X\mapsto X'$ such that $
\alpha^{x\lambda}=\alpha^{\lambda x^\varphi }$  for all $\alpha \in \Delta$ and  $x \in X$.

Let $G$ be a semiregular permutation group on a set $\Omega$ with $m$ orbits for some positive integer $m$. Set
\[ V=\cup_{i=1}^{m}G_{i} \mbox{ with } G_i=\{g_{i}\ :\ g\in G\}, \mbox{ and } R(G)=\{R(g)\ :\ g \in G \}, \]
where $R(g)$ is the right multiplication of $g$ on $V$, that is, $R(g): x_i \mapsto (xg)_i, \text{ for all } x  \in G \text{ and } 1\leq i \leq m$. Then $R(G)$ is a semiregular permutation group on $V$ with $m$ orbits: $G_1, G_2, \cdots, G_m$. For simplicity of notation, from now on and beyond, we write $R(G)$ instead of $R_m(G)$, as given in Section~\ref{Sec1}. The following lemma is standard, and for the integrity, we give a simple proof.

\begin{lem}
The permutation groups $G$ (on $\Omega$) and $R(G)$ (on $V$) are permutation isomorphic.
\end{lem}

\begin{proof}
Let $\Omega_1$, $\Omega_2$, $\cdots$, $\Omega_m$ be all orbits of $G$ on $\Omega$ and assume that $\a_i\in\Omega_i$ for each $1\leq i\leq m$. Since $G$ is semiregular, the map \begin{equation}\label{eq:lambda}
 \lambda : \alpha_{i}^{g}\mapsto g_{i},\ \forall g\in G,\ 1\leq i\leq m .
\end{equation}
is a bijection from $\Omega$ to $V$, which indeed induces a bijection from every $\Omega_i$ to $G_i$.
Let $\varphi$ be the isomorphism from $G$ to $R(G)$, defined by
\begin{equation}\label{eq:varphi}
\varphi: g\mapsto R(g), \forall g \in G.
\end{equation}
Note that  every element in $\Omega$ has the form $\a_i^h$ for some $1\leq i\leq m$ and $h\in G$.
Then for all $g\in G$ and $\a_i^h\in \Omega$, we have
\[(\alpha_i^h)^{g\lambda}=(\alpha_{i}^{hg})^{\lambda}=(hg)_{i}=(h_i)^{R(g)}=(\alpha_i^h)^{\lambda R(g)}=(\alpha_i^h)^{\lambda g^{\varphi}}.\]
This proves that $G$ (on $ \Omega$) and $R(G)$ (on $V$) are permutation isomorphic.
\end{proof}

Note that $\lambda$ and $\varphi$, given in Eqs~\eqref{eq:lambda} and~\eqref{eq:varphi}, induce the permutation isomorphic between $G$ and $R(G)$, that is, $\alpha^{g\lambda}=\alpha^{\lambda g^\varphi }$  for all $\alpha \in \Omega$ and  $g \in G$. It follows that $g=\lambda R(g)\lambda^{-1}$  for all $g \in G$, and then it is easy to check that  $$N_{\mathrm{S}_\Omega}(G)=\lambda N_{\mathrm{S}_V}(R(G))\lambda^{-1},$$ where $N_{\mathrm{S}_\Omega}(G)$ and $N_{\mathrm{S}_V}(R(G))$ are the normalizers of $G$ and $R(G)$ in $\S_\Omega$ and $\S_V$, respectively.

Thus, to determine $N_{\mathrm{S}_\Omega}(G)$, it suffices to determine the normalizer $N_{\mathrm{S}_V}(R(G))$ of $R(G)$ in $\mathrm{S}_V$. Denote by $C_{\mathrm{S}_V}(R(G))$ the centralizer of $R(G)$ in $\S_V$.

From now on, we consider $N_{\mathrm{S}_V}(R(G))$ and its subgroups. For simplicity, set
\[
N=N_{\mathrm{S}_V}(R(G)),\ C=C_{\mathrm{S}_V}(R(G)), \mbox { and }\ \G:=\{G_1,G_2,\ldots,G_m \}.\]
Then $\G$ is the orbit set of $R(G)$ on $V$. Since $R(G)$ is semiregular, $N$ is transitive on $V$, and since $R(G)\unlhd N$, $\G$ is a complete imprimitive block system of $N$. Let $K$ be the kernel of $N$ on $\G$. Then
$$K=\{ n\in N\ :\ G_i^n=G_i \mbox{ for every } 1\leq i\leq m\}.$$

For $g\in G$ and $1\leq i\leq m$, define $L_i(g)$ by
\[
 L_i(g):\ x_{i}\mapsto (g^{-1}x)_{i}, \text{ and } x_{j}\mapsto x_{j} \text{ for all } x\in G \text{ and } 1\leq j\leq m \text{ with } j\neq i.\]
Then $L_i(g)$ is a permutation on $V$. Set
 \[L_i=\{L_i(g): g\in G\} \text{ and } L=L_1L_2\cdots L_m.\]
Then $L_i$ is a permutation group on $V$, which is regular on $G_i$ and fixes $\bigcup\limits_{j\not=i}^{m}G_{j}$ pointwise. It is easy to see that for all $1\leq i,j\leq m$ with $i\neq j$, $L_i\cap L_j=1$ and $L_i$ and $L_j$ commutes pointwise. Thus, $$L = L_1\times L_2\times\cdots\times L_m.$$
Clearly, $|L_i|=|R(G)|=|G|$ and $|L|=|G|^m$.

Recall that $\S_m$ is the symmetric group on $\{1,2,\cdots,m\}$, and $\Aut(G)$ is the automorphism group of $G$.
We view them as permutation groups on $V$ as follows: for every $\sigma\in \S_{m}$, $\alpha \in \mathrm{Aut}(G)$, $i\in \{1,2,\ldots,m\}$ and $x\in G$:
\begin{align*}
&\sigma: x_i \mapsto x_{i^\sigma},\ \ \ \ \ \ \
\alpha: x_i \mapsto (x^\alpha)_i.\\
\end{align*}

\begin{theorem}\label{thC}
$C=L\rtimes \S_m$.
\end{theorem}

\begin{proof}
 For $g,h\in G$ and $1\leq i\le m$, we have $x_{i}^{L_i(h)R(g)}=(h^{-1}xg)_i=x_i^{R(g)L_i(h)}$, and  $x_{j}^{L_i(h)R(g)}=(xg)_j=x_j^{R(g)L_i(h)}$ for all $1\leq j\leq m$ with $j\not=i$. It follows that $L_i(h)R(g)=R(g)L_i(h)$, and hence $L=L_1\times L_2\times\cdots\times L_m\leq C$.

For $g\in G$ and $\s\in \S_{m}$, we have $x_{i}^{\s R(g)}=(xg)_{i^{\s}}=x_{i}^{R(g)\s}$ for all $1\leq i\leq m$ and $x\in G$. This implies that $\S_m\leq C$.

Note that $\S_m$ acts faithfully on $\G=\{G_1,G_2,\cdots,G_m\}$, and since $L$ fixes $G_i$ for all $1\leq i\leq m$, we have $L\cap \S_m=1$.

Take $g\
\in G$ and $\s\in \S_m$. For each $1\leq i\leq m$, we have $x_{j}^{\s^{-1}L_i(g)\s}=(g^{-1}x)_{j}=x_{j}^{L_{i^\s}(g)}$ when $j=i^\s$, and $x_{j}^{\s^{-1}L_i(g)\s}=x_{j}=x_{j}^{L_{i^\s}(g)}$ when $j\neq i^\s$. Then we have:
\begin{equation}\label{Eq1}
\s^{-1}L_i(g)\s=L_{i^\s}(g), \mbox{ for all } g\in G,\s\in \S_m \mbox{ and } 1\leq i\leq m.
\end{equation}
Thus, $(L_i)^\s=L_{i^\s}$, and hence $\S_m$ normalizes $L$. It follows $L\rtimes \S_{m}\leq C$.

Let $H=K\cap C$, that is, the kernel of $C$ on $\G=\{G_{1},\cdots, G_{m}\}$. Then $C/H\leq \S_m$.
Clearly, $L\leq H$. To prove $C=L\rtimes \S_{m}$, it suffices to prove $H=L $.

For every $1\leq i\leq m$, consider the constituent $H^{G_i}$ of $H$ on $G_i$, that is, the permutation group restricted by $H$ on $G_i$. Let $H_i$ be the kernel of $H$ acting on $G_i$. Then $H/H_i\cong H^{G_i}$. Clearly, $H_1\cap H_2\cap\cdots\cap H_m=1$. Since $R(G)$ is semiregular on $V$ with orbit set $\G$, the constituent $R(G)^{G_i}$ of $R(G)$ on $G_i$ is a regular permutation group on $G_i$. Note that $G_i$ is an orbit of $L_i$ and the constituent $L_i^{G_i}$ is also a regular permutation group on $G_i$, which centralizes $R(G)^{G_i}$.  By \cite[Lemma~4.2A]{Dixon}, $C_{\S_{G_i}}(R(G)^{G_i})=L_i^{G_i}$, and so $|C_{\S_{G_i}}(R(G)^{G_i})| = |L_i^{G_i}|=|G|$.

Since $H\leq C$, we have $H^{G_i}\leq C_{\S_{G_i}}(R(G)^{G_i})=L_i^{G_i}$ and hence $|H/H_i|=|H^{G_i}|\leq |G|$. Note that $H=H/(H_1\cap H_2\cap\cdots\cap H_m)$ is isomorphic to a subgroup of $H/H_1\times H/H_2\times\cdots\times H/H_m$. Thus, $|H|\leq |G|^m$.
Then it follows from $L\leq H$ and $\vert L \vert=|G|^m$ that $H=L$, as required.
\end{proof}

\medskip
By Theorem~\ref{thC}, the centralizer of a semiregular permutation group on $\Omega$ in $\S_{\Omega}$ is transitive (also see \cite[Exercise~4.2.5]{Dixon}). In particular, Theorem~\ref{thC} implies $C=L(G)$ for $m=1$, which is well-known (see \cite[Lemma~4.2A]{Dixon}).

\begin{theorem}\label{thN}
$N=(L_1\times\cdots\times L_m)\rtimes (\S_m\times \Aut(G))$.
\end{theorem}

\begin{proof}
Let $R(g)\in R(G)$ and $\alpha \in \Aut(G)$. Then we have
$$x_{i}^{\alpha^{-1}R(g)\alpha}=(x^{\alpha^{-1}}g)^{\alpha}_{i}=(xg^{\alpha})_{i}=x_{i}^{R(g^{\alpha})},\mbox{ for all } x\in G \mbox{ and } 1\leq i\leq m.$$
It follows
\begin{equation}\label{Eq2}
R(g)^{\alpha}=R(g^{\alpha}), \mbox{ for all } g\in G \mbox{ and } \alpha\in \Aut(G).
\end{equation}
This implies that $\Aut(G)$ normalizes $R(G)$, that is, $\Aut(G)\leq N$.

By Theorem~\ref{thC}, $C=L\rtimes \S_m\leq N$. Let $\s\in\S_m$ and $\a\in\Aut(G)$. For every $x\in G$ and $1\leq i\leq m$, we have $(x_i)^{\a\s}=(x^\a)_{i^\s}=(x_i)^{\s\a}$, and so $\S_m$ commutes with $\Aut(G)$ pointwise. Assume $\s=\a\in \S_m\cap \Aut(G)$. Then $x_{i^\s}=(x_i)^\s=(x_i)^\a=(x^\a)_i$, forcing $i^\s=i$ and $x^\a=x$. It follows $\s=\a=1$, and hence $\S_m\Aut(G)=\S_m\times\Aut(G)\leq N$.

Take $\a\in\Aut(G)$ and $L_i(g)\in L_i$ with $g\in G$. For $x\in G$, we have $(x_j)^{\a^{-1}L_i(g)\a}=(g^{-\a}x)_j=x_j^{L_i(g^\a)}$ for $j=i$, and $(x_j)^{\a^{-1}L_i(g)\a}=x_j=x_j^{L_i(g^\a)}$ for all $1\leq j\leq m$ with $j\not=i$. It follows
\begin{equation}\label{Eq3}
\a^{-1}L_i(g)\a=L_i(g^\a), \mbox{ for all } g\in G,\a\in \Aut(G) \mbox{ and } 1\leq i\leq m.
\end{equation}
This means $(L_i)^\a=L_i$, and hence $\Aut(G)$ normalizes every $L_i$. In particular, $\Aut(G)$ normalizes $L$, and so by Theorem~\ref{thC}, $\S_m\times \Aut(G)$ normalizes $L$.

Assume that $L_1(g_1)L_2(g_2)\cdots L_m(g_m)=\s\a$ for some $g_1,g_2,\cdots,g_m\in G$, $\s\in\S_m$ and $\a\in \Aut(G)$. For all $x\in G$ and $1\leq i\leq m$, we have $(g_i^{-1}x)_i=x_i^{L_1(g_1)L_2(g_2)\cdots L_m(g_m)}=x_i^{\s\a}=(x^\a)_{i^\s}$. This implies $i^\s=i$ and $g_i^{-1}x=x^\a$. Since $1\leq i\leq m$, we have $\s=1$, and by taking $x=1$ we have $g_i=1$. It follows $L_1(g_1)L_2(g_2)\cdots L_m(g_m)=1$, and hence $L\cap (\S_m\times \Aut(G))=1$. Thus, $L(\S_m\times \Aut(G))=L\rtimes (\S_m\times\Aut(G))\leq N$.

By the $N/C$ Theorem (see \cite{Dixon}), $|N|/|C|\leq |\Aut(R(G))|=|\Aut(G)|$. Since $C=L\rtimes \S_m$, we have $|N|\leq |C||\Aut(G)|=|L|\cdot |\S_m|\cdot |\Aut(G)|$, and since $L\rtimes (\S_m\times\Aut(G))\leq N$, we have $N=L\rtimes (\S_m\times\Aut(G))=(L_1\times L_2\times \cdots\times L_m)\rtimes (\S_m\times\Aut(G))$.
\end{proof}

\medskip
By Theorems~\ref{thC} and \ref{thN}, the group structures  of $C=(L_1\times L_2\times\cdots\times L_m)\rtimes \S_m$ and $N=(L_1\times\cdots\times L_m)\rtimes (\S_m\times \Aut(G))$ are completely determined by Eq~(\ref{Eq1}) and Eq~(\ref{Eq3}).

\medskip

By Eqs~(\ref{Eq2}) and (\ref{Eq3}), $\Aut(G)$ normalizes $R(G)$ and $L_i$ for every $1\leq i\leq m$. Let $1\leq r\leq m$. It is easy to see that $R(G)\cap (L_1\times\cdots\times L_{r-1}\times L_{r+1}\times\cdots\times L_m)=1$ by considering its action on $G_r$, and by Theorem~\ref{thC}, \[L_1\cdots L_{r-1}R(G)L_{r+1}\cdots L_m \Aut(G)=(L_1\times\cdots\times L_{r-1}\times R(G)\times L_{r+1}\times\cdots\times L_m)\rtimes \Aut(G).\]
 For $g\in G$, let $\Inn(g)$ be the inner automorphism of $G$ induced by $g$, that is, $x^{\Inn(g)}=x^g$ for all $x\in G$.
It is easy to check that $L_1(g)L_2(g)\cdots L_m(g)R(g)=\Inn(g)\in \Aut(G)$, and hence $(L_1\times L_2\times\cdots\times L_m)\rtimes \Aut(G)=(L_1\times\cdots\times L_{r-1}\times R(G)\times L_{r+1}\times\cdots\times L_m)\rtimes \Aut(G)$ for every $1\leq r\leq m$. Then Theorem~\ref{thN} implies the following.

\begin{cor}\label{c-R(G)} $N=((L_1\times\cdots \times L_{r-1}\times R(G)\times L_{r+1}\times\cdots\times L_m)\rtimes \Aut(G))\rtimes\S_m$ for every $1\leq r\leq m$.
\end{cor}

By Corollary~\ref{c-R(G)}, if $m=1$, then $N=R(G)\rtimes \Aut(G)$, which is the well-known holomorph of $G$ (see \cite[Corollary~4.2B]{Dixon}).
Recall that $K$ is the kernel of $N$ acting on $\G=\{G_1,G_2,\cdots,G_m\}$.

\begin{cor}\label{c-kernel}
The group  $K=L\rtimes\Aut(G)=(L_1\times\cdots \times L_{r-1}\times R(G)\times L_{r+1}\times\cdots\times L_m)\rtimes \Aut(G)$  for every $1\leq r\leq m$.
\end{cor}

For every $1\leq r\leq m$, Let $N_{1_r}$ and $K_{1_r}$ be the stabilizers of $1_r$ in $N$ and $K$ respectively, and $N_{G_r}$ the block stabilizer of $G_r$ in $N$.
The group $\S_m$ induces the symmetric group on $\mathcal{G}=\{G_1,G_2,\ldots,G_m\}$.
Let $(\S_{m})_{G_r}$ be the setwise stabilizer of $G_r$ in the group  $\S_m$.
By the definition of $\S_m$, $(\S_{m})_{G_r}\cong\S_{m-1}$ fixes each point in $G_r$.

\begin{cor}\label{c-stabilizer} For every $1\leq r\leq m$, we have
 $N_{1_r}=(L_1\times \cdots \times L_{r-1}\times L_{r+1}\times\cdots\times L_m)\rtimes ((\S_{m})_{G_r}\times \Aut(G))$, $N_{G_r}=(L_1\times \cdots\times L_m)\rtimes ((\S_{m})_{G_r}\times \Aut(G))$, and $K_{1_r}=(L_1\times \cdots \times L_{r-1}\times L_{r+1}\times\cdots\times L_m)\rtimes \Aut(G)$.
\end{cor}

\begin{proof}
By Eq~(\ref{Eq1}) and Eq~(\ref{Eq3}), $(\S_{m})_{G_r}$ and $\Aut(G)$ normalize $(L_1\times \cdots \times L_{r-1}\times L_{r+1}\times\cdots\times L_m)$, and by definitions of $L_i$ and $\Aut(G)$, both $(L_1\times \cdots \times L_{r-1}\times L_{r+1}\times\cdots\times L_m)$ and $\Aut(G)$ fix $1_r$. Since $(\S_{m})_{G_r}$ fixes $r$, it also fixes $1_r$, and by Theorem~\ref{thN}, $(L_1\times \cdots \times L_{r-1}\times L_{r+1}\times\cdots\times L_m)\rtimes ((\S_{m})_{G_r}\times \Aut(G))\leq N_{1_r}$. On the other hand, the transitivity of $N$ on $V$ implies $|N_{1_r}|=|N|/|V|=|G|^{m-1}(m-1)!|\Aut(G)|$, and hence $N_{1_r}=(L_1\times \cdots \times L_{r-1}\times L_{r+1}\times\cdots\times L_m)\rtimes ((\S_{m})_{G_r}\times \Aut(G))$.

Since $G_r$ is an imprimitive block of $N$, we have $N_{1_r}\leq N_{G_r}$. Clearly, $L_r\leq N_{G_r}$ is transitive on $G_r$, and hence $N_{G_r}=L_rN_{1_r}$. Since $|L_r|=|G|$, we have $N_{G_r}=(L_1\times \cdots\times L_m)\rtimes ((\S_{m})_{G_r}\times \Aut(G))$.

Since $L_r\leq K$, $K$ is transitive on $G_r$, and hence $K=L_rK_{1^r}$. By Corollary~\ref{c-kernel}, $K=L\rtimes\Aut(G)$, and since $L_r\cap K_{1^r}=1$, we have $|K_{1^r}|=|K|/|G|=|G|^{m-1}|\Aut(G)|$.
Note that $(L_1\times \cdots \times L_{r-1}\times L_{r+1}\times\cdots\times L_m)\rtimes \Aut(G)\leq K_{1_r}$. Since $|(L_1\times \cdots \times L_{r-1}\times L_{r+1}\times\cdots\times L_m)\rtimes \Aut(G)|=|G|^{m-1}|\Aut(G)|=|K_{1^r}|$, we have $K_{1_r}=(L_1\times \cdots \times L_{r-1}\times L_{r+1}\times\cdots\times L_m)\rtimes \Aut(G)$.
\end{proof}

\medskip
Let $m\geq 2$. For convenience, we write $V_m=V=G_1\cup G_{2}\cup\dots\cup G_{m}$ and $V_{m-1}=G_1\cup G_{2}\cup\dots\cup G_{m-1}$. Recall that $R(G)\leq \S_{V_m}$ is semiregular on $V_m$ with orbits set $\{G_1,G_2,\cdots,G_m\}$, and the constituent $R(G)^{V_{m-1}}$ of $R(G)$ on $V_{m-1}$ is a semiregular group on $V_{m-1}$ with orbits set  $\{G_1,G_2,\cdots,G_{m-1}\}$. We still use $R(G)$ to denote this constituent, that is $R(G)^{V_{m-1}}=R(G)$, and this makes no confusion in the context. By Theorem~\ref{thN},
$N_{S_{V_{m-1}}}(R(G))=(L_1\times\cdots\times L_{m-1})\rtimes (\S_{m-1}\times \Aut(G))$. By Corollary~\ref{c-stabilizer},  $N_{G_m}=(L_1\times \cdots\times L_m)\rtimes ((\S_{m})_{G_m}\times \Aut(G))$, and hence the kernel of $N_{G_m}$ on $V_{m-1}$ is $L_m$. It follows that
$|N_{G_m}^{V_{m-1}}|=|G|^{m-1}||\S_{m-1}||\Aut(G)|=|N_{S_{V_{m-1}}}(R(G))|$. Noting that $N_{G_m}=N_{V_{m-1}}$, we have the following result.

\begin{cor}\label{c-restriction}
$N_{S_{V_{m-1}}}(R(G))=N_{G_m}^{V_{m-1}}=N_{V_{m-1}}^{V_{m-1}}$.
\end{cor}

To end this section, we prove the following result for later use.

\begin{theorem}\label{SRGconj}
Let $G$ and $H$ are two isomorphic semiregular permutation groups on a finite set $\Omega$. Then $H$ and $G$ are conjugate in $\S_{\Omega}$.
\end{theorem}

\begin{proof}
 Since $G$ is semiregular, it has $ |\Omega|/|G|$ orbits on $\Omega$, and similarly, $H$ has ${|\Omega|}/{|H|}$ orbits. Since $G\cong H$, we may write    $m= {|\Omega|}/{|G|}= {|\Omega|}/{|H|}$.
Let $\overline{\Omega}=\{\Omega_1,\Omega_2,\cdots,\Omega_m\}$ and $\overline{\Delta}=\{\Delta_1,\Delta_2,
\cdots,\Delta_m\}$ be the orbit sets of $G$ and $H$ on $\Omega$, respectively. Then $\overline{\Omega}$ and $\overline{\Delta}$ are partitions of $\Omega$, that is, $\Omega=\Omega_1\cup\cdots\cup\Omega_m=\Delta_1\cup\cdots\cup\Delta_m$ with $\Omega_i\cap\Omega_j=\emptyset$ and $\Delta_i\cap\Delta_j=\emptyset$ for all $i\not=j$.
For every $1\leq i\leq m$, take $\a_i\in \Omega_i$ and $\b_i\in \Delta_i$. Then $\Omega_i=\a_i^G$ and $\Delta_i=\b_i^H$. Let $\a$ be an isomorphism from $G$ to $H$. Then $G^\a=H$. Define a map $\sigma$ on $\Omega$ as follows:
\begin{center}
$\sigma:\ \alpha_{i}^{g}\mapsto \beta_{i}^{g^{\alpha}}$, \ for all $g\in G \mbox{ and }\ 1\leq i \leq m$.
\end{center}
Let $\alpha_{i}^{g_{1}}=\alpha_{j}^{g_{2}}\in \Omega$ with $g_{1},g_{2}\in G$. Then $i=j$ because $\Omega_k=\a_k^G$ for all $1\leq k\leq m$ and $\overline{\Omega}$ is a partition of $\Omega$, and hence $g_1=g_2$ because $G$ is semiregular. It follows that $\beta_{i}^{g_1^{\alpha}}=\beta_{j}^{g_2^{\alpha}}$, and so $\sigma$ is well defined. On the other hand, let $\beta_{i}^{g_1^{\alpha}}=\beta_{j}^{g_2^{\alpha}}$. Since $H$ is semiregular, a similar argument above gives rise to $i=j$ and $g_1^\a=g_2^\a$, and since $\a$ is an isomorphism from $G$ to $H$, we have $g_1=g_2$. It follows that  $\alpha_{i}^{g_{1}}=\alpha_{j}^{g_{2}}$, that is, $\s$ is injective. Since $\Omega$ is finite, $\sigma$ is a bijection and hence $\s\in\S_{\Omega}$.

Recall that $G$ and $H$ are semiregular subgroups of $\S_{\Omega}$. Let $g\in G$. Take an arbitrary $\beta_{i}^{h}\in \Omega$ with $h\in H$. Since  $(\beta_i^h)^{\sigma^{-1}}=\alpha_i^{h^{\alpha^{-1}}}$, we have  $(\beta_{i}^{h})^{g^{\sigma}}=(\beta_{i}^{h})^{\sigma^{-1}g\sigma}=
(\alpha_{i}^{h^{\alpha^{-1}}})^{g\sigma}=\beta_{i}^{(h^{\alpha^{-1}}g)^{\alpha}}=
(\beta_{i}^{h})^{g^{\alpha}}$, and by the arbitrariness of $\beta_{i}^{h}$ in $\Omega$, we gave $g^{\sigma}=g^{\alpha}$. By the arbitrariness of $g$ in $G$, we have $G^{\sigma}=G^{\alpha}=H$, and hence $G$ is conjugate to $H$ in $\S_{\Omega}$.
\end{proof}

\section{$m$-Cayley digraph and its automorphisms \label{sec3}}

Let $m$ be a positive integer. For every $1\leq i,j \leq m$, let $S_{i,j}$ be subsets of $G$ with $1\notin S_{i,i}$. Let $V$ be the vertex set of the $m$-Cayley digraph $\Gamma=\Cay(G,S_{i,j}: 1\leq i,j\leq m)$ of $G$. Then
\begin{align*}
V=V(\Gamma)&=\bigcup_{1\leq i\leq m}  G_i, \text{ where } G_i=\{x_i :  x\in G\},\\
\mathrm{Arc} (\Gamma)&=\bigcup_{1\leq i,j\leq m} \{(x_i, (sx)_j) : s\in S_{i,j},x\in G\}.
\end{align*}
It is easy to see that a digraph $\Gamma$ is an $m$-Cayley digraph of $G$ if and only if $\Aut(\Gamma)$ has a subgroup isomorphic to $G$ acting semiregularly on the vertex set $V(\Gamma)$ with $m$ orbits (see \cite{DFS}).

The notations in Section~\ref{Sec2}, like $\G=\{G_1,G_2,\cdots,G_m\}$, $R(G)$, $N=N_{\S_V}(R(G))$, $C=C_{\S_V}(R(G))$, $L_i$, $\S_m$, $\Aut(G)$, $K$ (the kernel of $N$ on $\G$), and vertex stabilizer and block stabilizer of $N$ and $K$, will be used throughout the paper.

Now we describe the normalizer of $R(G)$ in $\Aut(\Gamma)$ and its subgroups. Write $R_1=\{g_1,g_2,\cdots,g_m\in G, \s\in\S_m\}$ and $R_2=\{g_1,g_2,\cdots,g_m\in G, \a\in \Aut(G)\}$. In this section we set

\begin{eqnarray*}
  {\tilde{N}}&=&\{L_1(g_1)\cdots L_m(g_m)\a\s\ |\ R_2,\s\in\S_m, S_{i^\s, j^\s}=(g _j^{-1}S_{i,j}g_i)^\a, 1\leq i,j\leq m\};\\
  {\tilde{C}}&=&\{L_1(g_1)\cdots L_m(g_m)\s\ |\ R_1, S_{i^\s, j^\s}=g_j^{-1}S_{i,j}g_i,1\leq i,j\leq m\};\\
      {\tilde{K}}& = &\{L_1(g_1)\cdots L_m(g_m)\a\ |\ R_2, S_{i,j}=(g _j^{-1}S_{i,j}g_i)^\a,1\leq i,j\leq m\};\\
   {\tilde{N}}_{1_r} &=&\{L_1(g_1)\cdots L_m(g_m)\a\s\ |\ R_2,g_r=1, \s\in(\S_{m})_{G_r}, S_{i^\s, j^\s}=(g _j^{-1}S_{i,j}g_i)^\a, 1\leq i,j\leq m\};\\
  {\tilde{N}}_{G_r}&=&\{L_{1}(g_{1})\cdots L_{m}(g_{m})\a\s\ |\ R_2, \s\in(\S_{m})_{G_r}, S_{i^\s, j^\s}=(g _j^{-1}S_{i,j}g_i)^\a, 1\leq i,j\leq m \};\\
   {\tilde{K}}_{1_r}& = & \{L_1(g_1)\cdots L_m(g_m)\a\ |\ R_2,g_r=1, S_{i, j}=(g _j^{-1}S_{i,j}g_i)^\a, 1\leq i,j\leq m\};\\
  \overline{\S}_m & = & \{\s\in \S_m\ |\ R_2, L_1(g_1)\cdots L_m(g_m)\a\s\in \tilde{N}\};\\
    \overline{\Aut}(G)& =&  \{\a\in \Aut(G)\ |\ R_1, L_1(g_1)\cdots L_m(g_m)\a\s\in \tilde{N}\}.
\end{eqnarray*}

\begin{theorem}\label{thNmC}
Let $\Gamma=\Cay(G,S_{i,j}: 1\leq i,j\leq m)$ and $A=\Aut(\Gamma)$. Then $N_A(R(G))=\tilde{N}$, $C_A(R(G))=\tilde{C}$, and the kernel of $\tilde{N}$ acting on $\{G_1,G_2,\cdots,G_m\}$ is $\tilde{K}$. As given above, $\tilde{K}_{1_r}$ and $\tilde{N}_{1_r}$, for every $1
\leq r\leq m$,  are stabilizers of $1_r$ in $\tilde{K}$ and $\tilde{N}$, respectively, and $\tilde{N}_{G_r}$ is the set-stabilizer of $G_r$ in  $\tilde{N}$. Furthermore, $\tilde{K}=R(G)\rtimes \tilde{K}_{1_r}$,  $\tilde{N}_{G_r}=R(G)\rtimes \tilde{N}_{1_r}$ and $\tilde{N}=\tilde{K}.\overline{\S}_m=\tilde{C}.\overline{\Aut}(G)$.
\end{theorem}
\begin{proof}
By Theorem~\ref{thN}, $N=N_{\S_{V(\Gamma)}}(R(G))=(L_1\times L_2\times\cdots\times L_m)\rtimes(\Aut(G)\times \S_m)$, the normalizer of $R(G)$ in $S_{V(\Gamma)}$. Clearly, $N_A(R(G))=N\cap A$.

Let $n\in N$. Then $n=L_1(g_1)\cdots L_m(g_m)\a\s$ for some $g_1,\cdots, g_m\in G,\a\in\Aut(G)$ and $\s\in\S_m$.
For any $(x_i,(s_{i,j}x)_j)\in \Arc(\Gamma)$, where $x\in G$, $1\leq i,j\leq m$ and $s_{i,j}\in S_{i,j}$, we have
\begin{equation}\label{EqninN}
(x_i,(s_{i,j}x)_j)^n=((g_i^{-1}x)^\a_{i^\s},(g_j^{-1}s_{i,j}x)^\a_{j^\s})=((g_i^{-1}x)^\a_{i^\s},((g_j^{-1}s_{i,j}g_i)^\a (g_i^{-1}x)^\a)_{j^\s}).
\end{equation}

If $n\in A$, by Eq~(\ref{EqninN}) we have $(g_j^{-1}s_{i,j}g_i)^\a \in  S_{i^\s, j^\s}$, and then $(g _j^{-1}S_{i,j}g_i)^\a\subseteq S_{i^\s, j^\s}$. Since $n$ maps $G_i$ to $G_{i^\s}$ and $G_j$ to $G_{j^\s}$, the number of out-neighbours of $1_i$ in $G_j$ is equal to the number of out-neighbours of $1_{i^\s}$ in $G_{j^\s}$, implying that $|S_{i,j}|=|S_{i^\s, j^\s}|$. Thus,  $|(g _j^{-1}S_{i,j}g_i)^\a|=|S_{i,j}|=|S_{i^\s, j^\s}|$, and therefore, $(g _j^{-1}S_{i,j}g_i)^\a=S_{i^\s, j^\s}$. On the other hand, if $(g _j^{-1}S_{i,j}g_i)^\a=S_{i^\s, j^\s}$ then Eq~(\ref{EqninN}) implies that $n$ maps an arc of $\Gamma$ to an arc, forcing $n\in A$. It follows that $n\in N\cap A$ if and only if  $(g _j^{-1}S_{i,j}g_i)^\a=S_{i^\s, j^\s}$, and by the definition of $\tilde{N}$, we have $N_A(R(G))=N\cap A=\tilde{N}$.

By Theorems~\ref{thC} and Corollary~\ref{c-kernel}, $C=(L_1\times L_2\times\cdots\times L_m)\rtimes \S_m$ and $K=(L_1\times L_2\times\cdots\times L_m)\rtimes \Aut(G)$ are the centralizer of $R(G)$ in $S_{V(\Gamma)}$ and the kernel of $N$ acting on the orbit set $\G=\{G_1,G_2,\cdots,G_m\}$ of $R(G)$, respectively. It follows that $C\cap A$ and $K\cap A$ are the centralizer of $R(G)$ in $A$ and the kernel of $A\cap N$ acting on $\G$, respectively. By taking $\a=1$ or $\sigma=1$ in Eq~(\ref{EqninN}), an argument similar to the above paragraph gives rise to $C_A(R(G))=C\cap A=\tilde{C}$, and $K\cap A=\tilde{K}$ as $A\cap N=\tilde{N}$.

Let $1\leq r\leq m$. By Corollary~\ref{c-stabilizer}, $N_{1_r}=(L_1\times \cdots \times L_{r-1}\times L_{r+1}\times\cdots\times L_m)\rtimes ((\S_{m})_{G_r}\times \Aut(G))$, $N_{G_r}=(L_1\times \cdots\times L_m)\rtimes ((\S_{m})_{G_r}\times \Aut(G))$, and $K_{1_r}=(L_1\times \cdots L_{r-1}\times L_{r+1}\times\cdots\times L_m)\rtimes \Aut(G)$. Then $N_{1_r}\cap A$ and $K_{1_r}\cap A$ are the stabilizers of $1_r$ in $N\cap A$ and $K\cap A$ respectively, and $N_{G_r}\cap A$ is the set stabilizer of $G_r$ in $N\cap A$.
Recall that  $N\cap A=\tilde{N}$ and $K\cap A=\tilde{K}$.
Using Eq~(\ref{EqninN}), an argument similar to the proof of $N_A(R(G))=\tilde{N}$ gives rise to
$\tilde{N}_{1_r}=N_{1_r}\cap A=\{L_1(g_1)\cdots L_m(g_m)\a\s\ |\ g_1,\cdots, g_m\in G,g_r=1, \a\in \Aut(G), \s\in(\S_{m})_{G_r}, S_{i^\s, j^\s}=(g _j^{-1}S_{i,j}g_i)^\a,1\leq i,j\leq m\}$,
${\tilde{K}}_{1_r}=K_{1^r}\cap A=\{L_1(g_1)\cdots L_m(g_m)\a\ |\ g_1,\cdots, g_m\in G,g_r=1, \a\in \Aut(G),
S_{i, j}=(g _j^{-1}S_{i,j}g_i)^\a,1\leq i,j\leq m\}$, and ${\tilde{N}}_{G_r}=N_{G_r}\cap A=\{L_{1}(g_{1})\cdots L_{m}(g_{m})\a\s\ |\ g_{1},\cdots g_{m}\in G, \a\in\Aut(G), \s\in(\S_{m})_{G_r}, S_{i^\s, j^\s}=(g _j^{-1}S_{i,j}g_i)^\a, 1\leq i,j\leq m \}$.

Since $R(G)\leq A$ fixes $\G$ pointwise, we have $R(G)\leq \tilde{K}$ and $R(G)\leq \tilde{N}_{G_r}$. Note that $R(G)$ is regular on $G_r$. By Frattini argument, $\tilde{K}=R(G)\rtimes \tilde{K}_{1_r}$ and  $\tilde{N}_{G_r}=R(G)\rtimes \tilde{N}_{1_r}$.

Recall that $N=(L_1\times L_2\times\cdots\times L_m)\rtimes(\Aut(G)\times \S_m)$. Take any two elements in $N$, say $L_1(g_1)\cdots L_m(g_m)\a\s$ and $L_1(g_1')\cdots L_m(g_m')\a'\s'$, where $g_1,\cdots, g_m,g_1',\cdots,g_m'\in G,\a,\a'\in\Aut(G)$ and $\s,\s'\in\S_m$. By Eqs~(\ref{Eq1}) and (\ref{Eq3}), $L_1(g_1)\cdots L_m(g_m)\a\s L_1(g_1')\cdots L_m(g_m')\a'\s'
=L_{1}(g_{1}(g_{1^{\s}}')^{\a^{-1}})\cdots L_{m}(g_{m}(g_{m^{\s}}')^{\a^{-1}})\a\a'\s\s'$. This implies that $f_A: L_1(g_1)\cdots L_m(g_m)\a\s\mapsto \a$ and $f_S: L_1(g_1)\cdots L_m(g_m)\a\s\mapsto \s$ are epimorphisms from $N$ to $\Aut(G)$ and from $N$ to $\S_m$, respectively. Clearly, $f_A$ and $f_S$ have kernels $\Ker(f_A)=C=(L_1\times L_2\times\cdots\times L_m)\rtimes \S_m$ and $\Ker(f_S)=K=(L_1\times L_2\times\cdots\times L_m)\rtimes \Aut(G)$.

The restriction of $f_A$ on $A\cap N$ is an epimorphism from $A\cap N=\tilde{N}$ to $\overline{\Aut}(G)$ with kernel $C\cap A=\tilde{C}$ because $\Ker(f_A)=C$. It follows that $\tilde{N}=\tilde{C}.\overline{\Aut}(G)$. Similarly, since $\Ker(f_S)=K$ we have $\tilde{N}=\tilde{K}.\overline{\S}_m$.
\end{proof}

\medskip
Let $m=1$. Then $\tilde{N}=\tilde{K}$. By Theorem~\ref{thNmC}, $\tilde{K}_{1_{1}}=\{\a\ |\ \a\in \Aut(G), S_{1,1}=S_{1,1}^\a\}=\Aut(G,S_{1,1})$ and $\tilde{N}=R(G)\rtimes \Aut(G,S_{1,1})$, which was given in~\cite{Godsil1}. Let $m=2$. Then $\tilde{N}=\tilde{K}$ or $\tilde{N}=\tilde{K}.\mz_2$. By Theorem~\ref{thNmC}, $\tilde{K}_{1_1}=\{L_2(g)\a\ |\ g\in G, \a\in \Aut(G), S_{1,1}=S_{1,1}^\a, S_{2,2}=(g^{-1}S_{2,2}g)^\a,
S_{1,2}=(g^{-1}S_{1,2})^\a, S_{2,1}=(S_{2,1}g)^\a\}=\{L_2(g)\a^{-1}\ |\ g\in G, \a\in \Aut(G), S_{1,1}^\a=S_{1,1}, S_{2,2}^\a=g^{-1}S_{2,2}g, S_{1,2}^\a=g^{-1}S_{1,2}, S_{2,1}^\a=S_{2,1}g\}$ and $\tilde{K}=R(G)\rtimes \tilde{K}_{1_1}$. Then \cite[Theorem 1.1]{ZF} follows from the fact $\tilde{N}=\tilde{K}$ or $\tilde{N}=\tilde{K}.\mz_2$. Let $m\geq 2$. Then \cite[Theorem 2.3]{HKM} follows from $N_A(R(G))=\tilde{N}$. Furthermore, \cite[Lemma 2.4]{HKM} follows from $\tilde{K}=R(G)\rtimes \tilde{K}_{1_{1}}$, \cite[Lemma 2.5]{HKM} follows from
$\tilde{N}_{1_1}=\{L_1(g_1)\cdots L_m(g_m)\a\s\ |\ g_1,\cdots, g_m\in G,g_1=1, \a\in \Aut(G), \s\in\S_{m-1}^1, S_{i^\s, j^\s}=(g _j^{-1}S_{i,j}g_i)^\a,1\leq i,j\leq m\}$, and \cite[Theorem 2.6]{HKM} follows from $\tilde{N}=\tilde{K}.\overline{\S}_m$ and $\tilde{N}_{G_1}=R(G)\rtimes \tilde{N}_{1_1}$.

Note that $N=K\rtimes\S_m$ is a semiproduct, but $\tilde{N}=\tilde{K}.\overline{\S}_m$ may not.
For example, let $p$ be a prime and $r\geq 2$ an integer, and let $m=p^r$ and $G=\langle a\rangle\cong\mz_p$. Take $S_{i,(i+1)}=\{1\}$ for $1\leq i\leq m-1$, $S_{m,1}=\{a\}$, and $S_{i,j}=\emptyset$ otherwise. Then  $\Gamma=\Cay(G,S_{i,j}: 1\leq i,j\leq m)$ is a directed cycle of length $p^{r+1}$. Thus,  $A=\Aut(\Gamma)=\mz_{p^{r+1}}$ and $\tilde{N}=A$. Clearly, the kernel $\tilde{K}$ of $\tilde{N}$ is $R(G)\cong\mz_p$. However, $\tilde{N}$ cannot be a semiproduct of $\tilde{K}=R(G)$ by a group because $\tilde{N}$ is cyclic.

\section{CI-properties related to $m$-Cayley  digraphs \label{sec4}}

Recall that a Cayley (di)graph $\mathrm{Cay}(G, S)$ is said to be CI if, for any
Cayley (di)graph $\mathrm{Cay}(G, T)\cong \mathrm{Cay}(G, S)$, there exists some  $\alpha \in \mathrm{Aut}(G)$ such that $T=S^{\alpha}$, and that a group $G$ is said to be  DCI  or  CI  if all Cayley  digraphs or Cayley graphs of $G$ are CI, respectively.
The following  problems concerning CI-properties of Cayley digraphs have received considerable attention in the literature:
\begin{itemize}
\item [] Problem (I): Which Cayley (di)graphs for a group $G$ are CI?
\item [] Problem (II): Which groups are DCI-groups or CI-groups?
\end{itemize}
In this section, we study Problem~(I) for $m$-Cayley  (di)graphs, and in the next section, we study Problem~(II) for $m$-Cayley (di)graphs.

We use the notations in Section~\ref{Sec2}, and  in this section, the normalizer $N$ of $R(G)$ in the symmetric group on $V=\cup_{i=1}^m G_i$ is important. By  Theorem~\ref{thN} and Corollary~\ref{c-kernel}, we have
$$N=(L_1\times\cdots\times L_m)\rtimes (\S_m\times \Aut(G)),$$
where $L_i$, $\mathrm{S}_m$, and $\mathrm{Aut}(G)$ are given in Section~\ref{Sec2}. To study problem~(I),  Babai's criterion is a useful tool, and in this section, we introduce some CI-properties of $m$-Cayley  digraphs and give analogues to Babai's criterion for each kind of CI-property of $m$-Cayley digraphs.
\subsection{mCI-(di)graphs\label{sec41}}
\begin{defi}\label{B.1}
An $m$-Cayley (di)graph $\Gamma$ is called $m$CI , if for every $m$-Cayley (di)graph $\Sigma$ of $G$ isomorphic to $\Gamma$, there is  $n\in N$ such that $\Gamma^n=\Sigma$.
\end{defi}

Recall that $V(\Gamma)=V=\cup_{i=1}^m G_i$. For $n\in N\leq \S_{V}$, the (di)graph $\Gamma^n$ has vertex set $V(\Gamma^n)=\{v^n: v \in V(\Gamma)\}=\{v^n: v \in V\}=V(\Gamma)=V$ and arc set $\mathrm{Arc}(\Gamma^n)=\{(u^n,v^n) : \ (u,v)\in \mathrm{Arc}(\Gamma)\}$.

Let $n\in N$ and $1\leq r\leq m$. By Corollary~\ref{c-R(G)}, there exist $g,g_1,g_2,\cdots, g_m\in G$ with $g_r=1$, $\a\in \Aut(G)$, and $\s\in\S_m$, such that $$n=L_1(g_1)\cdots L_{r-1}(g_{r-1})R(g)L_{r+1}(g_{r+1})\cdots L_m(g_m)\a\s.$$ Note that $R(g)L_i(g_i)=L_i(g_i)R(g)$ and $\Gamma^{R(g)}=\Gamma$ as $R(g)\in R(G)\leq \Aut(\Gamma)$. From Eq~(\ref{EqninN}), it is easy to see that
\begin{equation}\label{Eqgamman}
\Gamma^n=\Cay(G,L_{i^\s,j^\s}: 1\leq i,j\leq m), \mbox{ where } L_{i^\s,j^\s}=(g_j^{-1}S_{i,j}g_i)^\a.
\end{equation}

\medskip

From the above Equation we immediately have the following proposition.

\begin{prop}\label{mCImPCI}  Let $1\leq r\leq m$. Then an $m$-Cayley (di)graph $\Cay(G,S_{i,j}: 1\leq i,j\leq m)$ is $m$CI if and only if whenever $\Cay(G,S_{i,j}: 1\leq i,j\leq m)$ is isomorphic to an $m$-Cayley (di)graph  $\Cay(G,T_{i,j}: 1\leq i,j\leq m)$, there exist $\alpha\in \Aut(G), \sigma\in \S_{m}, g_1,\ldots,g_m\in G$  with $g_{r}=1$, such that $T_{i^{\sigma},j^{\sigma}}=g_jS^{\alpha}_{i,j}g_{i}^{-1}$.
\end{prop}

By Proposition~\ref{mCImPCI}, a $1$-Cayley (di)graph $\Cay(G,S_{1,1})$ is $1$CI if and only if, for any $\Cay(G,T_{1,1})\cong \Cay(G,S_{1,1})$, there exists some $\alpha \in \mathrm{Aut}(G)$ such that $T_{11}=S_{11}^{\alpha}$.
Thus, the $1$CI is exactly the CI for Cayley (di)graph.

Again by Proposition~\ref{mCImPCI}, a $2$-Cayley (di)graph $\Cay(G,S_{i,j}:1\leq i,j\leq 2)$ is $2$CI if and only if, for any $\Cay(G,T_{i,j}:1\leq i,j\leq 2)$ isomorphic to $\Cay(G,S_{i,j}:1\leq i,j\leq 2)$, there exists some $\alpha \in \mathrm{Aut}(G)$, $\sigma\in \mathrm{S}_m$, and $g\in G$ such that
either $T_{1,1}=S_{1,1}^\a$, $T_{2,2}=g^{-1}S_{2,2}^\a g$, $T_{1,2}=g^{-1}S_{1,2}^\a$, $T_{2,1}=S_{2,1}^\a g$, or $T_{2,2}=S_{1,1}^\a $, $T_{1,1}=g^{-1}S_{2,2}^\a g$, $T_{2,1}=g^{-1}S_{1,2}^\a$ and $T_{1,2}=S_{2,1}^\a g$.
If $\Cay(G,S_{i,j}:1\leq i,j\leq 2)$ is a graph, then $T_{1,2}=g^{-1}S_{1,2}^\a$ if and only if $T_{2,1}=S_{2,1}^\a g$. Thus, the $2$CI is exactly the SCI for $2$-Cayley graph, as defined in Arezoomand and Taeri~\cite{Arez}.

\medskip
The following is a generalization of the well-know Babai criterion of Cayley (di)graph to $m$-Cayley (di)graph.

\begin{theorem}\label{BabaiSimilar} Let $\Gamma=\Cay(G,S_{i,j}: 1\leq i,j\leq m)$ be an
$m$-Cayley (di)graph. Then
 $\Gamma$ is $m$CI if and only if every semiregular group of $\Aut(\Gamma)$ isomorphic to $G$ is conjugate to $R(G)$ in $\Aut(\Gamma)$.
\end{theorem}

\demo Recall that $V(\Gamma)=\bigcup_{i=1}^m G_i$ and $N=N_{\S_{V(\Gamma)}}(R(G))=(L_1\times L_2\times\cdots\times L_m)\rtimes(\Aut(G)\times \S_m)$. Set $A=\Aut(\Gamma)$.

To prove the necessity, assume that $\Gamma$ is $m$CI and  that $H\cong R(G)$ is an $m$-semiregular subgroup of automorphisms of $\Gamma$. By Lemma~\ref{SRGconj}, $R(G)=H^\s$ for some $\s\in\S_{V(\Gamma)}$.
Note that $\s$ is an isomorphism from $\Gamma$ to $\Gamma^\s$. Since $H$ is a semiregular subgroup of $\Aut(\Gamma)$, $H^\s$ is a semiregular subgroup of $\Aut(\Gamma^\s)$, and since $R(G)=H^\s$, we have $R(G)\leq \Aut(\Gamma^\s)$. This implies that $\Gamma^\s$ is an $m$-Cayley (di)graph of $G$. Since $\Gamma$ is $m$CI, there is $n\in N$ such that $\Gamma^n=\Gamma^\s$. Thus, $n\s^{-1}\in\Aut(\Gamma)=A$ and $R(G)^{n\s^{-1}}=R(G)^{\s^{-1}}=H$, that is, $H$ and $R(G)$ are conjugate in $A$.

To prove the sufficiency, assume that every semiregular group of automorphisms of $\Gamma$ isomorphic to $G$ is conjugate to $R(G)$ in $\Aut(\Gamma)$. Let $\Sigma=\Cay(G,T_{i,j}: 1\leq i,j\leq m)$ be an $m$-Cayley (di)graph of $G$, and let $\s$ be an isomorphism from $\Gamma$ to $\Sigma$. Then $\Gamma^\s=\Sigma$, and since $R(G)\leq \Aut(\Sigma)$ is semiregular on $V(\Sigma)$, $R(G)^{\s^{-1}}$ is a semiregular subgroup of $A$, which is isomorphic $R(G)$. By assumption, there is $\a\in\Aut(\Gamma)$ such that $R(G)^\a=R(G)^{\s^{-1}}$. It follows that $\a\s\in N$ and $\Gamma^{\a\s}=\Gamma^\s=\Sigma$, that is, $\Gamma$ is $m$CI.
\hfill\qed

\medskip
For $m=1$, Theorem~\ref{BabaiSimilar} implies the well-known Babai criterion for Cayley (di)graph (see \cite[Lemma 3.1]{Babai2}), and for $2$-Cayley graph, Theorems~\ref{BabaiSimilar} and \ref{SRGconj} imply \cite[Theorem~A]{Arez}.

\subsection{mPCI-(di)graphs\label{sec42}}
Let $\Gamma$ and $\Sigma$ be two $m$-PCayley (di)graphs of a group $G$. An isomorphism $\sigma$ from $\Gamma$ to $\Sigma$ is called a {\em $p$-isomorphism} if $\{G_1,G_2,\cdots,G_m\}$ is $\sigma$-invariant, that is, $\sigma$ maps every $G_i$ for $1\leq i\leq m$ to some $G_j$, and $\Gamma$ and $\Sigma$ are called $p$-isomorphic if there is a $p$-isomorphic from $\Gamma$ and $\Sigma$, denoted by $\Gamma\cong_{p}\Sigma$. Now we state the definition of $m$PCI for $m$-PCayley (di)graph, which is also motivated from \cite{Arez,Ko1,Ko2}.

\begin{defi}\label{PCI}
An $m$-PCayley (di)graph $\Gamma$ of a group $G$ is called $m$PCI, if for any $m$-PCayley (di)graph $\Sigma$ $p$-isomorphic to $\Gamma$, there is $n\in N$ such that $\Gamma^n=\Sigma$.
\end{defi}

It is worth mentioning that Arezoomand~\cite[Lemma 4.4]{Arez} proved that if $\Gamma$ and $\Sigma$ are isomorphic $2$-PCayley graphs on a group, then $\Gamma\cong_{p}\Sigma$. Since an m-PCayley (di)graph is an m-Cayley (di)graph, Proposition~\ref{mCImPCI} implies that $2$PCI-graph is exactly $0$-type SCI-graph, as defined in \cite[Definition 4.1(i)]{Arez}.

\medskip
For $m$-PCayley (di)graphs, similar to Theorem~\ref{BabaiSimilar} we have the following generalization of Babai criterion of Cayley (di)graph to $m$-PCayley (di)graph. Recall that $\G=\{G_1,G_2,\cdots,G_m\}$.

\begin{theorem}\label{PCI-BabaiSimilar}
Let $\Gamma=\Cay(G,S_{i,j}: 1\leq i,j\leq m)$ be an
$m$-PCayley (di)graph. Then
$\Gamma$ is $m$PCI if and only if every semiregular group of $\Aut(\Gamma)$ isomorphic to $G$, with the same orbit set $\G$ as $R(G)$, is conjugate to $R(G)$ in $\Aut(\Gamma)$.
\end{theorem}

\demo  Write $N=N_{\S_{V(\Gamma)}}(R(G))$ and  $A=\Aut(\Gamma)$.

To prove the necessity, assume that $\Gamma$ is $m$PCI and $H$ is an $m$-semiregular subgroup of automorphisms of $\Gamma$ with orbit set $\G$. By Lemma~\ref{SRGconj}, $R(G)=H^\s$ for some $\s\in\S_{V(\Gamma)}$, and since $H$ has orbit set $\G$, $H^\s$ has orbit set $\{G_1^\s,G_2^\s,\cdots,G_m^\s\}$, forcing that $\G=\{G_1^\s,G_2^\s,\cdots,G_m^\s\}$, that is, $\s$ is a $p$-isomorphism from $\Gamma$ to $\Gamma^\s$. This also implies that $\Gamma^\s$ is an $m$-PCayley (di)graph of $G$, and since $\Gamma$ is $m$PCI, there is $n\in N$ such that $\Gamma^n=\Gamma^\s$. It follows $n\s^{-1}\in A$ and  $R(G)^{n\s^{-1}}=R(G)^{\s^{-1}}=H$.

To prove the sufficiency, assume that every semiregular group of automorphisms of $\Gamma$  isomorphic to $G$, with orbit set $\G$, is conjugate to $R(G)$ in $A$. Let $\Sigma=\Cay(G,T_{i,j}: 1\leq i,j\leq m)$ be an $m$-PCayley (di)graph of $G$ such that $\Gamma\cong_{p}\Sigma$. Then there is a $p$-isomorphism from $\Gamma$ to $\Sigma$, say $\s$, and so $\Gamma^\s=\Sigma$ and $\G$ is $\s$-invariant. It follows that $R(G)$ and $R(G)^{\s^{-1}}$ are two isomorphic semiregular subgroups of $A$ with the same orbit set $\G$, and by assumption, there is $\a\in A$ such that $R(G)^\a=R(G)^{\s^{-1}}$. Then $\a\s\in N$ and $\Gamma^{\a\s}=\Gamma^\s=\Sigma$.
 \hfill\qed

\medskip

To end this section, we prove that $m$PCI properties of an $m$-PCayley (di)graph and its multipartite complement are equivalent, which is important for classifying $m$(P)(D)CI-groups in Section~\ref{sec5}. Let $\Gamma=\Cay(G,S_{i,j}: 1\leq i,j\leq m)$ be an $m$-PCayley (di)graph. Then  $S_{i,i}=\emptyset$ for $1\leq i\leq m$. Set $T_{i,i}=\emptyset$ and $T_{i,j}=G \backslash{S_{i,j}}$ for all $1\leq i,j\leq m$ with $i\not=j$. Then $\Cay(G,T_{i,j}: 1\leq i,j\leq m)$ is called {\em the multipartite complement} of $\Gamma$, denoted by $\Gamma^{mc}$. It is easy to see that $(\Gamma^{mc})^{mc}=\Gamma$ and $\Gamma^{mc}$ is an $m$-PCayley (di)graph. Recall that one of $\Gamma$ and the component $\Gamma^c$ of $\Gamma$ must be connected. However, it is possible that both  $\Gamma$  and $\Gamma^{mc}$ are disconnected.

\begin{lem}\label{cmofdigraph}
Let $m\geq 2$ be an integer and let $\Gamma=\Cay(G,S_{i,j}: 1\leq i,j\leq m)$ and $\Sigma=\Cay(G,T_{i,j}: 1\leq i,j\leq m)$ be $m$-PCayley digraphs.
\begin{itemize}
\item[\rm (1)] For $n\in N=N_{\S_{V(\Gamma)}}(R(G))$, we have $\Gamma^n=\Sigma$ if and only if $(\Gamma^{mc})^n=\Sigma^{mc}$.
\item[\rm (2)] If both $\Gamma$ and $\Gamma^{mc}$ are not connected, then $m=2$, $|G|$ is even and $\Gamma\cong 2K_{|G|/2,|G|/2}$.
\end{itemize}
\end{lem}

\begin{proof}
Let $n\in N=N_{\S_{V(\Gamma)}}(R(G))$ and $\Gamma^n=\Sigma$. Since $(\Gamma^{mc})^{mc}=\Gamma$, to prove (1), it suffices to prove that  $(\Gamma^{mc})^n=\Sigma^{mc}$.
By Theorem~\ref{thN}, $n=L_1(h_1)\cdots L_m(h_m)\a\s$ for some $h_1,\cdots, h_m\in G,\a\in\Aut(G)$ and $\s\in\S_m$.  Let $1\leq i,j\leq m$ with $i\not=j$. By Proposition~\ref{mCImPCI},  $\Gamma^n=\Sigma$ implies that $T_{i^\s,j^\s}=(h_j^{-1}S_{i,j}h_i)^\a$, and since $G^\a=G$, we have $G\backslash{T_{i^\s,j^\s}}=(h_j^{-1}(G\backslash{S_{i,j}})h_i)^\a$. Again by Proposition~\ref{mCImPCI}, $(\Gamma^{mc})^n=\Sigma^{mc}$, as required.

\smallskip
Now we prove (2). Assume that $\Gamma$ has $t$ components, say $\Gamma_1, \Gamma_2,\cdots,\Gamma_t$, where $t\geq 2$. Then all $V(\Gamma_j)$, $1\leq j\leq t$, are blocks of $\Aut(\Gamma)$. Set $G_i^j=G_i\cap V(\Gamma_j)$ for $1\leq i\leq m$ and $1\leq j\leq t$.
Then for each $1\leq i\leq m$, $G_i=\cup_{j=1}^t G_i^j$ is a disjoint union.

\medskip
\noindent {\bf Claim:}  $G_i^j\not=\emptyset$ for every $1\leq i\leq m$ and $1\leq j\leq t$, and each $G_i$ cannot be contained in a component of $\Gamma$ or $\Gamma^{mc}$.

Suppose to the contrary that $G_i^j=\emptyset$ for some given $i$ and $j$ with $1\leq i\leq m$ and $1\leq j\leq t$. Then we may let $G_i^k\not=\emptyset$ for every $1\leq k\leq \ell$, and $G_i^k=\emptyset$ for every $\ell+1\leq k\leq t$, where $1\leq \ell<t$. It follows that $G_i=\cup_{j=1}^\ell G_i^j\subseteq \cup_{j=1}^\ell V(\Gamma_j)$. Since $R(G)\leq \Aut(\Gamma)$ is transitive on every $G_i$ and $V(\Gamma_j)$ is a block of $\Aut(\Gamma)$ for every $1\leq j\leq t$, $R(G)$ fixes $\cup_{j=1}^\ell V(\Gamma_j)$ setwise, and hence $\cup_{j=1}^\ell V(\Gamma_j)$ is a union of some $G_k$ with $1\leq k\leq m$. Since $\ell<t$, $\Gamma^{mc}$ is connected because every vertex in $\cup_{j=1}^\ell V(\Gamma_j)$ is adjacent to every vertex in $V(\Gamma)\backslash \cup_{j=1}^\ell V(\Gamma_j)$ in $\Gamma^{mc}$, a contradiction. Thus, $G_i^j\not=\emptyset$ for every $1\leq i\leq m$ and $1\leq j\leq t$, and
$G_i=\cup_{j=1}^t G_i^j$ is a disjoint union.

Since $R(G)\leq \Aut(\Gamma)$, we have $|G_i^1|=|G_i^2|=\cdots=|G_i^t|$ for every $1\leq i\leq m$, and since $|G_i|=|G|$, all $|G_i^j|$ are equal one another. Suppose that $G_i\subseteq V(\Gamma_j)$ for a given $j$ with $1\leq j\leq t$. Then $G_i=G_i^j$, and since $|G_i^1|=|G_i^2|=\cdots=|G_i^t|$ and
$G_i=\cup_{j=1}^t G_i^j$ is a disjoint union, we have $t=1$, contradicting $t\geq 2$. Thus, $G_i$ cannot be contained in a component of $\Gamma$. Since $(\Gamma^{mc})^{mc}=\Gamma$ and $R(G)\leq \Aut(\Gamma^{mc})$, a similar argument above implies that $G_i$ cannot be contained in a component of $\Gamma^{mc}$. This establishes the claim.

\medskip
Suppose $t\geq 3$. For any $3\leq j\leq t$, $G_1^1$ and $G_1^j$ are in the same component of $\Gamma^{mc}$ because every vertex in $G_1^1$ is adjacent to every vertex in $G_2^2$ and every vertex in $G_2^2$ is adjacent to every vertex in $G_1^j$, and $G_1^1$ and $G_1^2$ are in the same component of $\Gamma^{mc}$ because every vertex in $G_1^1$ is adjacent to every vertex in $G_2^j$ and every vertex in $G_2^j$ is adjacent to every vertex in $G_1^2$. It follows that $G_1$ is contained in one component of $\Gamma^{mc}$, contradicting Claim. Thus, $t=2$ and $\Gamma=\Gamma_1\cup \Gamma_2$.

Suppose $m\geq 3$. Then $G_1^1$ and $G_1^2$ are contained in the same component of $\Gamma^{mc}$ because every vertex in $G_1^1$ is adjacent to every vertex in $G_2^2$, every vertex in $G_2^2$ is adjacent to $G_m^1$, and every vertex in $G_m^1$ is adjacent to every vertex in $G_1^2$. This implies that $G_1$ is contained in one component of $\Gamma^{mc}$ as $t=2$, contradicting Claim. Thus, $m=2$ and both $\Gamma_1$ and $\Gamma_2$ are bipartite. Furthermore, it is easy to see that if $\Gamma_1$ is not a complete bipartite graph, then $\Gamma^{mc}$ is connected, a contradiction. This yields that both $\Gamma_1$ and $\Gamma_2$ are complete bipartite graphs, of which each is isomorphic to $K_{|G|/2,|G|/2}$. This completes the proof.
\end{proof}

From Lemma~\ref{cmofdigraph}~(1), we have the following corollary.

\begin{cor}\label{B.3-C}
Let $m\geq 2$ be integer and let $\Gamma=\Cay(G,S_{i,j}: 1\leq i,j\leq m)$ be an $m$-PCayley digraph. Then $\Gamma$ is mPCI if and only if $\Gamma^{mc}$  is mPCI.
\end{cor}

\section{CI- and DCI-groups related to $m$-Cayley digraphs\label{sec5}}

Similar to DCI-groups and CI-groups for Cayley (di)graphs, we have the following definitions.

\begin{defi}\label{B.3}
A finite group $G$ is called $m$DCI (resp. $m$CI), if every $m$-Cayley digraph (resp. $m$-Cayley graph) of $G$ is $m$CI. Furthermore, $G$ is called $m$PDCI (resp. $m$PCI), if every $m$-PCayley digraph (resp. $m$-PCayley graph) of $G$ is $m$PCI.
\end{defi}

In this section, our main purpose is to classify $m$DCI- and $m$CI-groups for each $m\geq 2$, and $m$PDCI- and $m$PCI-groups for each  $m\geq 4$. Recall that a graph is viewed as a digraph, and an $m$-PCayley (di)graph is an $m$-Cayley (di)graph, which imply the following proposition.

\begin{prop}\label{relation}
Let $G$ be a finite group. Then
\begin{enumerate}
 \item[\rm (1)] If $G$ is $m$DCI or $m$PDCI then $G$ is $m$CI or $m$PCI, respectively;
  \item[\rm (2)] If $G$ is $m$DCI or $m$CI then $G$ is $m$PDCI or $m$PCI, respectively.
\end{enumerate}
\end{prop}

The following theorem is a reduction theorem, a key result to classify $m$(P)(D)CI-groups.

\begin{theorem}\label{mtom-1}
Let $m\geq 2$ be an integer and let $G$ be a finite group.
If $G$ is $m$DCI, $m$CI,  $m$PDCI or $m$PCI, then $G$ is $(m-1)$DCI, $(m-1)$CI, $(m-1)$PDCI or $(m-1)$PCI, respectively.
\end{theorem}

\begin{proof}
Write $V_m=G_1\cup G_{2}\cup\dots\cup G_{m}$ and $V_{m-1}=G_1\cup G_{2}\cup\dots\cup G_{m-1}$.

Assume that $G$ is $m$DCI. Let $\Gamma=\Cay(G,S_{i,j}: 1\leq i,j\leq m-1)$ and $\Sigma=\Cay(G,T_{i,j}: 1\leq i,j\leq m-1)$ be two isomorphic $(m-1)$-Cayley digraphs of $G$.
To prove that $G$ is $(m-1)$DCI, we only need to prove that there exists some $n\in N_{V_{m-1}}(R(G))$ such that $\Gamma^n=\Sigma$. By Eq.~(\ref{Eqgamman}), $\Gamma^n=\Sigma$ if and only if $(\Gamma^c)^n=\Sigma^c$, and hence we may assume that $\Gamma$ and $\Sigma$ are connected.

Note that $S_{i,j}$ and $T_{i,j}$ are given for every $1\leq i,j\leq m-1$. We extend them to $1\leq i,j\leq m$ by setting $S_{\ell,m}=S_{m,\ell}=S_{m,m}=T_{\ell,m}=T_{m,\ell}=T_{m,m}=\emptyset$ for each
$1\leq \ell\leq m-1$. Let $\Gamma_1=\Cay(G,S_{i,j}: 1\leq i,j\leq m)$ and $\Sigma_1=\Cay(G,T_{i,j}: 1\leq i,j\leq m)$. Then $\Gamma_1=\Gamma\cup G_m$, $\Sigma_1= \Sigma\cup G_m$, and so $\Gamma_1$ and $\Sigma_1$ are isomorphic $m$-Cayley digraphs. Set $N=N_{V_m}(R(G))$. Since $G$ is an $m$DCI-group, there is $n_1\in N$ such that $\Gamma_1^{n_1}=\Sigma_1$. Since $\Sigma$ and $\Gamma$ are connected, $n_1$ fixes $G_m$ setwise, that is, $n_1\in N_{G_m}=N_{V_{m-1}}$, and hence $\Gamma^{n_1}=\Sigma$. By Corollary~\ref{c-restriction},
$N_{V_{m-1}}(R(G))=N_{G_m}^{V_{m-1}}=N_{V_{m-1}}^{V_{m-1}}$. Let $n$ be the restriction of $n_1$ on $V_{m-1}$. We have $n\in N_{V_{m-1}}(R(G))$ and $\Gamma^n=\Sigma$, as required.

Assume that $G$ is $m$CI. Then a similar argument to $m$DCI implies that $G$ is $(m-1)$CI.

Now assume that $G$ is $m$PDCI. To prove that $G$ is $(m-1)$PDCI, let $\Gamma=\Cay(G,S_{i,j}: 1\leq i,j\leq m-1)$ and $\Sigma=\Cay(G,T_{i,j}: 1\leq i,j\leq m-1)$ be two $p$-isomorphic $(m-1)$-PCayley digraphs of $G$. It suffices to show that there is $n\in N_{V_{m-1}}(R(G))$ such that $\Gamma^n=\Sigma$. By Lemma~\ref{cmofdigraph}~(2), we may assume that $\Gamma$ is either connected, or $\Gamma\cong 2K_{|G|/2,|G|/2}$ where $m=2$ and $|G|$ is even. Then the similar argument to the above $m$DCI implies that $G$ is $(m-1)$PDCI. Furthermore, if $G$ is $m$PCI, then a similar argument to $m$PDCI implies that $G$ is $(m-1)$PCI.
\end{proof}

\subsection{$m$CI- and $m$DCI-groups}

In this subsection, we classify $m$CI-, $m$DCI-groups for every $m\geq 2$. By definition, $1$DCI-groups and $1$CI-groups, which are the well-known DCI-groups and CI-groups for Cayley (di)graphs.
The classifications of DCI-groups and CI-groups are long-standing open problems.
Therefore, we assume that $m\geq 2$.

\begin{theorem}\label{mCI-groups} Let $m\geq 2$ be an integer and let $G$ be a finite group. Then
  \begin{itemize}
    \item[\rm (1)] $G$ is $m$CI if and only if either $m=2$ and $G=1$ or $\mz_3$, or $m\geq 3$ and $G=1$;
    \item[\rm (2)]  $G$ is $m$DCI if and only if $G=1$.
  \end{itemize}
\end{theorem}
\begin{proof}
We first prove part~(1). For $m=2$, it follows from \cite[Theorem B]{Arez} that $G$ is $2$CI if and only if $G=1$ or $\mz_3$. Let $m\geq 3$. If $G$ is $m$CI, then, by Theorem~\ref{mtom-1}, $G$ is $2$CI, implying that either $G=1$ or $G=\mz_3$. On the other hand, Theorem~\ref{BabaiSimilar} implies that $G=1$ is an $m$CI-group.

Suppose $G=\mz_3=\langle x\rangle$. For all $1\leq i,j\leq 3$ with $i\not=j$, take  $S_{i,i}=\emptyset$ and $S_{i,j}=\{1\}$, and $T_{i,i}=\{x,x^{2}\}$ and $T_{i,j}=\emptyset$. Then $\Gamma=\Cay(G,S_{i,j}: 1\leq i,j\leq 3)$ and $\Sigma=\Cay(G,T_{i,j}: 1\leq i,j\leq 3)$ are union of three cycles of length $3$, and hence $\Gamma\cong\Sigma$. Since $|S_{i,i}|\neq|T_{j,j}|$ for all $1\leq i,j\leq 3$, Proposition~\ref{mCImPCI} implies that $\mathbb{Z}_{3}$ cannot be $3$CI, and this, together with Theorem~\ref{mtom-1}, implies that $\mathbb{Z}_{3}$ cannot be $m$CI for any $m\geq 3$. This completes the proof of part~(1).

Assume that $G$ is $m$DCI. By Proposition~\ref{relation}, $G$ is $m$CI, and by part~(1), either $m=2$ and $G=1$ or $\mz_3$, or $m\geq 3$ and $G=1$. On the other hand,  Theorem~\ref{BabaiSimilar} implies that $G=1$ is an $m$DCI-group for every integer $m\geq 1$. Thus, we may let $m=2$ and $G=\mz_3=\langle x\rangle$. Then $V(\Gamma)=G_1\cup G_2$. Take  $S_{1,1}=S_{2,2}=\{x\}$, $S_{1,2}=\{1,x,x^2\}$ and $S_{2,1}=\emptyset$. Let $\Gamma=\Cay(G,S_{i,j}: 1\leq i,j\leq 2)$. Then every vertex in $G_1$ has in-valency $1$ and out-valency $4$, and  every vertex in $G_2$ has in-valency $4$ and out-valency $1$. Each of the induced subgraph $[G_1]$ and $[G_2]$ is a directed cycle of length $3$. Thus, $\Aut(\Gamma)=\langle a\rangle\times\langle b\rangle\cong\mz_3\times\mz_3$, where $a=(1_1,x_1,x^2_1)$ and $b=(1_2,x_2,x^2_2)$. Clearly, $\langle ab\rangle$ and $\langle ab^{-1}\rangle$ are two semiregular subgroups of $\Aut(\Gamma)$, which are not conjugate in $\Aut(\Gamma)$. By Theorem~\ref{BabaiSimilar}, $\Gamma$ is not $2$CI, and hence $\mz_3$ is not a $2$DCI-group. This completes the proof of part~(2).
\end{proof}

\medskip

\subsection{$m$PCI- and $m$PDCI-groups}

In this subsection, we study $m$PCI- and $m$PDCI-groups.
Since $1$-PCayley digraph of a group $G$ is the empty graph with $|G|$ vertices, all finite groups are  $1$PDCI- and $1$PCI-groups. Thus, we assume $m\geq 2$.  By Proposition~\ref{relation}, $m$PDCI-groups are $m$PCI-groups, and by Theorem~\ref{mtom-1}, the main work is to classify $m$PCI-groups with $m$ as small as possible. In this subsection, we  classify $m$PDCI- and $m$PCI-groups for each  $m\geq 4$.

It is well known that every subgroup of a CI-group is CI (see \cite[Lemma 3.2]{Babai}).
We generalize this result to $m$PCI-group for $m\geq 3$.

\begin{theorem}\label{subgroupsmPCI}
Let $m\geq 3$ be an integer. Then every subgroup of a finite $m$PCI-group is $m$PCI.
\end{theorem}

\begin{proof}
Let $G$ be a finite $m$PCI-group and let $H\leq G$. Let $\Gamma=\Cay(H,S_{i,j}: 1\leq i,j\leq m)$ be a $m$-PCayley graph of $H$.
Then $V(\Gamma)=\cup_{i=1}^mH_i$ with $H_i=\{h_i\ |\ h\in H\}$. To prove that $\Gamma$ is $m$PCI, by Lemma~\ref{cmofdigraph} and Corollary~\ref{B.3-C} we may assume that $\Gamma$ is connected.

To avoid confusion of symbols, we write $R_{\Gamma}(H)$ as the right multiplication of $H$ on $V(\Gamma)$ (similarly, $R_{\Gamma}(h)$ is the right multiplication of $h \in H$).
Set $\H=\{H_1,H_2,\cdots,H_m\}$, the orbit set of $R_{\Gamma}(H)$, and denote by $\Aut(\Gamma)_{(\H)}$ the subgroup of $\Aut(\Gamma)$ fixing $\H$ pointwise. Let $L$ be a semiregular automorphism group of $\Gamma$ which is isomorphic to $R_{\Gamma}(H)$ and has the same orbit set $\H$. Then $R_{\Gamma}(H), L\leq \Aut(\Gamma)_{(\H)}$. To finish the proof, by Theorem~\ref{PCI-BabaiSimilar} we only need to show that $L$ and $R_{\Gamma}(H)$ are conjugate in $\Aut(\Gamma)$.

Let $\Sigma=\Cay(G,S_{i,j}: 1\leq i,j\leq m)$ and $t=|G:H|$. Let $g_1=1,g_2,
\cdots,g_t$ be a right transversal of $H$ in $G$, that is, $G=\cup_{k=1}^tHg_k$ with $Hg_i\not=Hg_j$ for all $i\not=j$. Then $V(\Sigma)=\cup_{i=1}^mG_i$ with $G_i=\cup_{k=1}^t(Hg_k)_i$, where $(Hg_k)_i=\{(hg_k)_i\ |\ h\in H\}$. Since $\Gamma$ is connected, it is a component of $\Sigma$, and $\Sigma$ has $t$ components: $\Sigma_1=\Gamma,\Sigma_2=\Gamma^{R(g_1)},\cdots,\Sigma_t=\Gamma^{R(g_t)}$, where  $$V(\Sigma_k)=\cup_{i=1}^m(Hg_k)_i \mbox{ and } A(\Sigma_k)=\{((hg_k)_i,(h'g_k)_j)\ |\ h,h'\in H, (h_i,h'_j)\in A(\Gamma),1
\leq i,j\leq m\}.$$  Then $\Sigma$ is a disjoint union of $\Sigma_1,\Sigma_2,\cdots, \Sigma_t$.
Let $$C\Sigma=\{V(\Sigma_1),V(\Sigma_2),\cdots,V(\Sigma_t)\}.$$
Then $C\Sigma$ is a complete imprimitive block system of $\Aut(\Sigma)$ because $\Gamma$ is connected.  Recall that  $\G=\{G_1,G_2,\cdots,G_m\}$ is the orbit set of $R(G)$ on $V(\Sigma)$, and $\Aut(\Sigma)_{(\G)}$ is the subgroup of $\Aut(\Sigma)$ fixing $\G$ pointwise.

It is well known that $\Aut(\Sigma)_{(\G)}$ is the wreath product $\Aut(\Gamma)_{(\H)}\wr \S_t$ of $\Aut(\Gamma)_{(\H)}$ by $\S_t$, as explained below for details. The symmetric group $\S_t$ on $\{1,2,\cdots,t\}$ is viewed as a subgroup of $\Aut(\Sigma)$ by, for every $\s\in\S_t$,
\[
\s:  (hg_k)_i \mapsto (hg_{k^\s})_i \text{ for all } 1\leq i\leq m, 1\leq k\leq t \text{ and } h\in H.
\]
Then $\S_t\leq \Aut(\Sigma)_{(\G)}$ and $\S_t$ induces the symmetric group on  $\{V(\Sigma_1),V(\Sigma_2),\cdots,V(\Sigma_t)\}$, which implies that $\S_t$ acts on $C\Sigma$ faithfully. Clearly, $R(G)\leq \Aut(\Sigma)_{(\G)}$.

Now we extend $\Aut(\Gamma)_{(\H)}$ as a permutation group on $V(\Sigma)$ such that $\Aut(\Gamma)_{(\H)}$ fixes $V(\Sigma_k)$ pointwise for all $2\leq k\leq t$. Then $\Aut(\Gamma)_{(\H)}\leq \Aut(\Sigma)_{(\G)}$, and hence $\Aut(\Gamma)_{(\H)}^{R(g_k)}\leq \Aut(\Sigma)_{(\G)}$, where $\Aut(\Gamma)_{(\H)}^{R(g_k)}$ fixes $V(\Sigma_i)$ pointwise for every $i\not=k$ as $R(g_k)$ maps $\Sigma_1$ to $\Sigma_k$.

Let $\s\in\S_t$. Since $R(g_k)\s$ and $R(g_{k^\s})$ map $\Sigma_1$ to $\Sigma_{k^\s}$, both $\Aut(\Gamma)_{(\H)}^{R(g_k)\s}$ and $\Aut(\Gamma)_{(\H)}^{R(g_{k^\s})}$ fix $V(\Sigma_i)$ pointwise for all $i\not=k^\s$. For $\a^{R(g_k)\s}\in \Aut(\Gamma)_{(\H)}^{R(g_k)\s}$ with $\a\in\Aut(\Gamma)_{(\H)}$ and for
$(hg_{k^\s})_i\in V(\Sigma_{k^\s})$ with $h\in H$ and $1\leq i\leq m$, let $h_i^\a=h(i)_i$ for some $h(i)\in H$. Then  $(hg_{k^\s})_i^{\a^{R(g_k)\s}}=(hg_{k^\s})_i^{\s^{-1}R(g_k^{-1})\a R(g_k)\s}=(h(i)g_{k^\s})_i=(hg_{k^\s})_i^{\a^{R(g_{k^\s})}}$. It follows that
$$(\a^{R(g_k)})^\s=\a^{R(g_{k^{\s}})}, \mbox{ for all }\a\in \Aut(\Gamma)_{(\H)}, 1\leq k\leq t \mbox{ and } \s\in\S_t.$$
This implies $(\Aut(\Gamma)_{(\H)}^{R(g_k)})^\s=\Aut(\Gamma)_{(\H)}^{R(g_{k^\s})}$ for every $\s\in\S_t$ and $1\leq k\leq t$, and then we have
\begin{equation*}
\Aut(\Gamma)_{(\H)}\wr \S_t:=(\Aut(\Gamma)_{(\H)}\times \Aut(\Gamma)_{(\H)}^{R(g_2)}\times \cdots\times \Aut(\Gamma)_{(\H)}^{R(g_t)})\rtimes \S_t= \Aut(\Sigma)_{(\G)},
\end{equation*}
where the last equality follows from the fact that the restriction of $\Aut(\Sigma)_{(\G)}$ on $\Sigma_1$ is $\Aut(\Gamma)_{(\H)}$ and $\Aut(\Gamma)_{(\H)}\times \Aut(\Gamma)_{(\H)}^{R(g_2)}\times \cdots\times \Aut(\Gamma)_{(\H)}^{R(g_t)}$ is the kernel of $\Aut(\Sigma)_{(\G)}$ acting on $C\Sigma$.

For $X\leq \Aut(\Gamma)_{(\H)}$ and $T\leq\S_t$, we also have the wreath product
\begin{equation}\label{XwrT}
X\wr T:=(X^{R(g_1)}\times X^{R(g_2)}\times \cdots\times X^{R(g_t)})\rtimes T,
\end{equation}
where $(\a^{R(g_k)})^\s=\a^{R(g_{k^\s})}$ for all $\a\in X$ and $\s\in T$. Furthermore, $X\times X^{R(g_2)}\times \cdots\times X^{R(g_t)}$ is the kernel of $X\wr T$ acting on $C\Sigma$.

Let $M$ be the subgroup of $\S_t$ induced by $R(G)$ on $C\Sigma$. By the transitivity of $R(G)$, $M$ is transitive on $C\Sigma$. Note that $R_{\Gamma}(H)\leq \Aut(\Gamma)_{(\H)}$ fixes $V(\Sigma_k)$ pointwise for all $2\leq k\leq t$. Since $g,g_2g,\cdots,g_tg$ is also a right transversal of $H$ in $G$ for any $g\in G$, we conclude that $R(G)$ normalizes $R_{\Gamma}(H)\times R_{\Gamma}(H)^{R(g_2)}\times \cdots\times R_{\Gamma}(H)^{R(g_t)}$.
Hence, $R(G)(R_{\Gamma}(H)\times R_{\Gamma}(H)^{R(g_2)}\times \cdots\times R_{\Gamma}(H)^{R(g_t)})\leq \Aut(\Sigma)_{(\G)}$, which induces the subgroup $M$ of $\S_t$ on   $C\Sigma$. Since $V(\Sigma_k)=\cup_{i=1}^m(Hg_k)_i$, the kernel of $R(G)$ acting on $C\Sigma$ is as the same as the kernel of $R(G)$ acting on the set of right cosets of $H$ in $G$, which is $R(K)$ with $K$ as the largest normal subgroup of $G$ contained in $H$. Then $R(K)\leq R_{\Gamma}(H)\times R_{\Gamma}(H)^{R(g_2)}\times \cdots\times R_{\Gamma}(H)^{R(g_t)}$ and hence we may write $R(G)(R_{\Gamma}(H)\times R_{\Gamma}(H)^{R(g_2)}\times \cdots\times R_{\Gamma}(H)^{R(g_t)})/(R_{\Gamma}(H)\times R_{\Gamma}(H)^{R(g_2)}\times \cdots\times R_{\Gamma}(H)^{R(g_t)})=M$. By Eq~(\ref{XwrT}), $(R_{\Gamma}(H)\times R_{\Gamma}(H)^{R(g_2)}\times \cdots\times R_{\Gamma}(H)^{R(g_t)})M/(R_{\Gamma}(H)\times R_{\Gamma}(H)^{R(g_2)}\times \cdots\times R_{\Gamma}(H)^{R(g_t)})=M$. It follows that
$R(G)(R_{\Gamma}(H)\times R_{\Gamma}(H)^{R(g_2)}\times \cdots\times R_{\Gamma}(H)^{R(g_t)})=(R_{\Gamma}(H)\times R_{\Gamma}(H)^{R(g_2)}\times \cdots\times R_{\Gamma}(H)^{R(g_t)})M$, forcing $R(G)\leq (R_{\Gamma}(H)\times R_{\Gamma}(H)^{R(g_2)}\times \cdots\times R_{\Gamma}(H)^{R(g_t)})M=R_{\Gamma}(H)\wr M$.

Recall that $L\leq \Aut(\Gamma)_{(\H)}$ is semiregular and isomorphic to $R_{\Gamma}(H)$ with the orbit set $\H$. By Eq~(\ref{XwrT}), we have the wreath product $L\wr M=(L\times L^{R(g_2)}\times \cdots\times L^{R(g_t)})\rtimes M$. Let $\b$ be an isomorphism from $R_{\Gamma}(H)$ to $L$. Then $\b$ can be extended to an isomorphism from $R_{\Gamma}(H)\wr M$ to $L\wr M$ by defining $m^\b=m$ and $(R_{\Gamma}(h)^{R(g_k)})^\b=(R_{\Gamma}(h)^\b)^{R(g_k)}$ for all $h\in H$, $m\in M$ and $1\leq k\leq t$, and for convenience, this extended automorphism is still denoted by $\b$.

Set $I=R(G)^\b$. Then $I\leq (R_{\Gamma}(H)\wr M)^\beta=L\wr M$ as $R(G)\leq R_{\Gamma}(H)\wr M$.
Since $L\leq  \Aut(\Gamma)_{(\H)}$ and $\Aut(\Gamma)_{(\H)}\wr \S_t=\Aut(\Sigma)_{(\G)} \leq \mathrm{Aut}(\Sigma)$, we see that $L\leq \mathrm{Aut}(\Sigma)$.
Now we claim that $I$ is semiregular and has the same orbit set $\G$ as $R(G)$.

For any given $\Sigma_i$ and $\Sigma_j$ with $1\leq i,j\leq t$, there is $R(g)\in R(G)$ such that $V(\Sigma_i)^{R(g)}=V(\Sigma_j)$. Since  $R(G)\leq R_{\Gamma}(H)\wr M$, we have $R(g)=xm$ for some $x\in R_{\Gamma}(H)\times R_{\Gamma}(H)^{R(g_2)}\times \cdots\times R_{\Gamma}(H)^{R(g_t)}$ and $m\in M$. Then $V(\Sigma_i)^{m}=V(\Sigma_i)^{xm}=V(\Sigma_j)$, and since $R(g)^\b=x^\b m$ with $x^\b\in L\times L^{R(g_2)}\times \cdots\times L^{R(g_t)}$, we have $V(\Sigma_i)^{R(g)^\b}=V(\Sigma_i)^{x^\b m}=V(\Sigma_i)^{m}=V(\Sigma_j)$. Thus, $I$ is transitive on  $C\Sigma$.

Since $R(G)_{V(\Sigma_1)}^{V(\Sigma_1)}=R_{\Gamma}(H)$ and $R_{\Gamma}(H)^\b=L$, we have $I_{V(\Sigma_1)}^{V(\Sigma_1)}=L$, which implies that $I_{V(\Sigma_1)}$ has $m$ orbits on $V(\Sigma_1)$, that is, $\H=\{H_1,H_2,\cdots,H_m\}$. Since $I$ is transitive on  $C\Sigma$, every orbit of $I$ has length $|G|$ as $|I|=|G|=|R(G)^\b|$. Thus, $I$ is semiregular, and since $L\wr M$ has orbit set $\G$, $I$ has the orbit set $\G$, as claimed.

Since $G$ is $m$PCI, there is $\a\in\Aut(\Sigma)$ such that $R(G)=I^\a$. Since $C\Sigma$ is a complete imprimitive block system of $\Aut(\Sigma)$, $V(\Sigma_1)^\a=V(\Sigma_i)$ for some $1\leq i\leq t$. Since $R(G)$ is transitive on $C\Sigma$, there is $g\in G$ such that $V(\Sigma_1)^{R(g)}=V(\Sigma_i)$, and since $\Sigma_1=\Gamma$, we have  $V(\Gamma)^{\a R(g^{-1})}=V(\Gamma)$, that is, $\a R(g^{-1})\in \Aut(\Sigma)$ fixes the component $\Gamma$. Write $\g=\a R(g^{-1})|_{V(\Gamma)}$, the restriction of $\a R(g^{-1})$ on $V(\Gamma)$.
Then $\g\in\Aut(\Gamma)$.

Since $\a R(g^{-1})$ fixes the component $\Gamma$, we get  $(I^{\a R(g^{-1})})_{V(\Gamma)}=(I_{V(\Gamma)})^{\a R(g^{-1})}$. Since $R(G)=I^\a$, we have $I^{\a R(g^{-1})}=R(G)$. It follows that $R(G)_{V(\Gamma)}=(I_{V(\Gamma)})^{\a R(g^{-1})}$ and hence $R_{\Gamma}(H)=R(G)_{V(\Gamma)}^{V(\Gamma)}=(I_{V(\Gamma)}^{\a R(g^{-1})})^{V(\Gamma)}=(I_{V(\Gamma)}^{V(\Gamma)})^{\g}=L^\g$, that is, $R_{\Gamma}(H)$ and $L$ are conjugate in $\Aut(\Gamma)$, as required.
\end{proof}

The proof in Theorem~\ref{subgroupsmPCI} does not work for $m=2$, because one cannot assume that $\Gamma$ is always connected by Lemma~\ref{cmofdigraph}.

\medskip
The following result shows that the elements of an $m$PCI-group for $m\geq 4$ are very restrictive.

\begin{lem}\label{pci2346}
Let $m\geq 4$ be a positive integer and $G$ a finite group. If $G$ is $m$PCI, then every element of $G$ has order $1,2,3,4$ or $6$.
\end{lem}
\begin{proof}
Let $G$ be a finite $m$PCI-group with $m\geq4$. Suppose on the contradiction that there is an element $x\in G$ with $o(x)=5$ or $\geq 7$, where $o(x)$ is the order of $x$ in $G$.
Let $\phi(k)$ be the number of elements coprime to $k$ in $\mz_k$. Then $\phi(o(x))\geq 3$, and hence $\langle x\rangle\cong\mz_k$ has $\phi(k)$ generators, that is $x^\ell$ for all $\ell\in\mz_k^*$, where $\mz_k^*$ is the multiplicative group of numbers coprime to $k$ in $\mz_k$.

Set $H=\langle x\rangle$. Take $S_{1,2}=S_{4,3}=\{1,x\}$, $S_{2,1}=S_{3,4}=\{1,x^{-1}\}$, $S_{2,3}=H=S_{3,2}$, and $S_{i,j}=\emptyset$ for all other $1\leq i,j\leq 4$. Let $\Gamma=\Cay(H,S_{i,j}: 1\leq i,j\leq 4)$. Then $V(\Gamma)=\cup_{i=1}^4H_i$ with $H_i=\{x_i\ |\ x\in H\}$. Clearly, $\Gamma=[H_1\cup H_2]\cup [H_3\cup H_4]\cup [H_2\cup H_3]$, where  $[H_1\cup H_2]=(1_1,x_2,x_1,x^2_2,x^2_1,\cdots,x^{k-1}_2,x^{k-1}_1,1_2)$ and $[H_3\cup H_4]=(1_4,x_3,x_4,x^3_2,x^2_4,\cdots,x^{k-1}_3,x^{k-1}_4,1_3)$ are $2k$-cycles, and $[H_2\cup H_3]\cong K_{k,k}$ is the complete bipartite graph of order $2k$ with partite sets $H_2$ and $H_3$, see below Figure~\ref{n4-pcay}. Recall that $\Aut(\Gamma)_{(\H)}$ is the subgroup of $\Aut(\Gamma)$ fixing $\H=\{H_1,H_2,H_3,H_4\}$ pointwise.

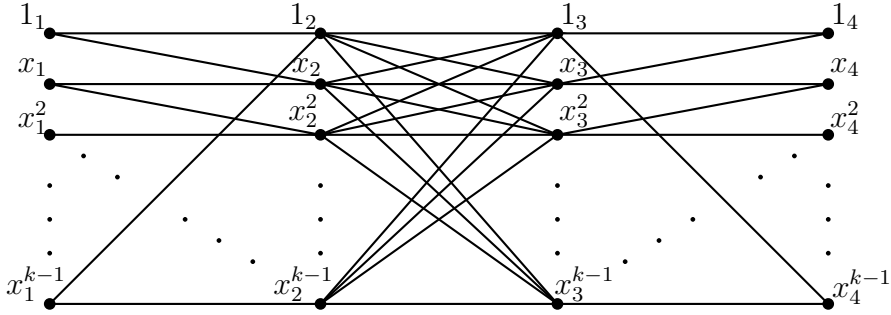
\begin{figure}[htp!]
  \centering
\begin{tikzpicture}[scale=0.45]
\draw [thick] (15,4)--(7,4)--(15,2.5)--(7,2.5)--(15,1)--(7,1);
\draw [thick] (15,-4)--(7,-4)--(15,4);
\draw [thick] (22,4)--(30,4)--(22,2.5)--(30,2.5)--(22,1)--(30,1);
\draw [thick] (22,4)--(30,-4)--(22,-4);

\fill (29,0.375) circle (2pt);
\fill (28,-0.25) circle (2pt);
\fill (26,-1.5) circle (2pt);
\fill (25,-2.125) circle (2pt);
\fill (24,-2.75) circle (2pt);

\fill (8,0.375) circle (2pt);
\fill (9,-0.25) circle (2pt);
\fill (11,-1.5) circle (2pt);
\fill (12,-2.125) circle (2pt);
\fill (13,-2.75) circle (2pt);

\filldraw[thick,fill=black] (7,4) circle (4pt);
\node at (6.5,4.5) {$1_1$};
\filldraw[thick,fill=black] (7,2.5) circle (4pt);
\node at (6.5,3) {$x_1$};
\filldraw[thick,fill=black] (7,1) circle (4pt);
\node at (6.5,1.5) {$x_1^2$};
\fill (7,-0.5) circle (2pt);
\fill (7,-1.5) circle (2pt);
\fill (7,-2.5) circle (2pt);
\filldraw[thick,fill=black] (7,-4) circle (4pt);
\node at (6.6,-3.4) {$x_1^{k-1}$};

\filldraw[thick,fill=black] (15,4) circle (4pt);
\node at (14.5,4.5) {$1_2$};
\filldraw[thick,fill=black] (15,2.5) circle (4pt);
\node at (14.6,3) {$x_2$};
\filldraw[thick,fill=black] (15,1) circle (4pt);
\node at (14.5,1.7) {$x_2^2$};
\fill (15,-0.5) circle (2pt);
\fill (15,-1.5) circle (2pt);
\fill (15,-2.5) circle (2pt);

\filldraw[thick,fill=black] (15,-4) circle (4pt);
\node at (14.5,-3.4) {$x_2^{k-1}$};

\filldraw[thick,fill=black] (22,4) circle (4pt);
\node at (22.5,4.5) {$1_3$};
\filldraw[thick,fill=black] (22,2.5) circle (4pt);
\node at (22.5,3) {$x_3$};
\filldraw[thick,fill=black] (22,1) circle (4pt);
\node at (22.5,1.7) {$x_3^2$};
\fill (22,-0.5) circle (2pt);
\fill (22,-1.5) circle (2pt);
\fill (22,-2.5) circle (2pt);
\filldraw[thick,fill=black] (22,-4) circle (4pt);
\node at (22.8,-3.4) {$x_3^{k-1}$};

\filldraw[thick,fill=black] (30,4) circle (4pt);
\node at (30.5,4.5) {$1_4$};
\filldraw[thick,fill=black] (30,2.5) circle (4pt);
\node at (30.5,3) {$x_4$};
\filldraw[thick,fill=black] (30,1) circle (4pt);
\node at (30.5,1.5) {$x_4^2$};
\fill (30,-0.5) circle (2pt);
\fill (30,-1.5) circle (2pt);
\fill (30,-2.5) circle (2pt);
\filldraw[thick,fill=black] (30,-4) circle (4pt);
\node at (31,-3.5) {$x_4^{k-1}$};

\draw [thick] (15,4)--(22,4)--(15,2.5)--(22,2.5)--(15,1)--(22,1)--(15,4);
\draw [thick]
(15,4)--(22,2.5);
\draw [thick]
(15,2.5)--(22,1);
\draw [thick]
(15,1)--(22,4);
\draw [thick]
(15,-4)--(22,-4);
\draw [thick]
(15,4)--(22,-4);
\draw [thick]
(15,2.5)--(22,-4);
\draw [thick]
(15,1)--(22,-4);
\draw [thick]
(22,4)--(15,-4);
\draw [thick]
(22,2.5)--(15,-4);
\draw [thick]
(22,1)--(15,-4);
\end{tikzpicture}
\caption{$4$-PCayley graph of $\mathbb{Z}_k=\langle x\rangle$}
\label{n4-pcay}
\end{figure}

Then $\Aut(\Gamma)_{(\H)}=\langle \a_1,\b_1\rangle \times \langle \a_2,\b_2\rangle\cong D_{2k}\times D_{2k}$ and $A:=\Aut(\Gamma)=\Aut(\Gamma)_{(\H)}\rtimes\langle \g\rangle$, where
\begin{align*}
&\begin{cases}
&\alpha_1=(1_1,x_1,\cdots,x_1^{k-1})(1_2,x_2,\cdots,x_2^{k-1}), \ \
\alpha_2=(1_3,x_3,\cdots,x_3^{k-1})(1_4,x_4,\cdots,x_4^{k-1}),\\
&\gamma=(1_2,1_3)(x_2,x_3)\cdots (x_2^{k-1},x_3^{k-1})(1_1,1_4)(x_1,x_4)\cdots (x_1^{k-1},x_4^{k-1}),\\
\end{cases}
\end{align*}
and for every $i\in \mz_k$, $\b_1$ interchanges $x^i_1$ and $x^{-i-1}_1$ and also $x^i_2$ and $x^{-i}_2$, and $\b_2$ interchanges $x^i_4$ and $x^{-i-1}_4$ and also $x^i_3$ and $x^{-i}_3$.
Note that $\b_1$ fixes $1_2$ and $\b_2$ fixes $1_3$.

Clearly, $\langle \a_1\a_2^\ell\rangle$ for all $\ell\in\mz_k^*$ are distinct semiregular subgroups of $A$ isomorphic to $\mz_k$, of which all have the orbit set $\H$. Let $L=\langle \a_1\a_2\rangle$. Then $\a_1,\a_2,\b_1\b_2,\g\in N_A(L)$, the normalizer of $L$ in $A$. Thus, $|N_A(L)|\geq 4k^2$, and since $|A|=8k^2$, we have $|A:N_A(L)|\leq 2$, which implies that the conjugate class of $L$ in $A$ contains at most $2$ subgroups. Since  $\phi(k)\geq 3$, all semiregular subgroups of $A$ isomorphic to $\mz_k$, with orbit set $\H$, cannot be conjugate in $A$. By Theorem~\ref{PCI-BabaiSimilar}, $\Gamma$ is not $4$PCI, and hence $H$ is not $4$PCI.
By Theorem~\ref{mtom-1}, $H$ cannot be $m$PCI for all $m\geq 4$, and by Theorem~\ref{subgroupsmPCI}, $G$ is not $m$PCI, contrary to assumption. It follows that $\phi(k)\leq 2$ for all $x\in G$ with $o(x)=k$, and hence $x$ has order $1,2,3,4$ or $6$.
\end{proof}

\medskip
Now we deal with some small non-$m$PCI-groups.

\begin{lem}\label{non-PCIs}
Each of $\mz_2^2$, $\mz_4$, $\mz_6$ and $\mz_3^2$ is not $4$PCI, and each of $\mz_3$ and $D_6$ is not $6$PCI.
\end{lem}

\begin{proof}
Assume $G=\langle x\rangle\times \langle y\rangle\cong \mathbb{Z}_2\times \mathbb{Z}_2$. Take $S_{1,2}=S_{2,1}=S_{3,4}=S_{4,3}=\{1,x\}$, $S_{1,3}=S_{3,1}=S_{2,4}=S_{4,2}=\{1,y\}$, $S_{2,3}=S_{3,2}=G$, and $S_{i,j}=\emptyset$  for all other $1\leq i,j\leq 4$. Let $\Gamma=\Cay(G,S_{i,j}: 1\leq i,j\leq 4)$. Then $V(\Gamma)=\cup_{i=1}^4G_i$ with $G_i=\{x_i\ |\ x\in G\}$, and $\G=\{G_1,G_2,G_3,G_4\}$. Clearly, $\Gamma=[G_1\cup G_2]\cup[G_1\cup G_3]\cup [G_2\cup G_3]\cup[G_2\cup G_4]\cup [G_3\cup G_4]$, where the induced subgraphs $[G_1\cup G_2]$, $[G_1\cup G_3]$, $[G_2\cup G_4]$ and $[G_3\cup G_4]$ are union of two $4$-cycles and $[G_2\cup G_3]$ is the complete bipartite graph of order $8$. It is easy to see that $|\Aut(\Gamma):\Aut(\Gamma)_{(\G)}|=4$, $\Aut(\Gamma)_{(\G)}\cong \mz_2^4$ and  $\alpha=(1_1,xy_1)(x_1,y_1)(1_2,xy_2)(x_2,y_2)(1_3,xy_3)(x_3,y_3)(1_4,xy_4)(x_4,y_4)$ lies in the center of $\Aut(\Gamma)_{(\G)}$, that is, $\a\in Z(\Aut(\Gamma)_{(\G)})$. Let
\begin{align*} \beta_1=(1_1,y_1)(x_1,xy_1)(1_2,y_2)(x_2,xy_2), & \ \  \beta_4=(1_1,xy_1)(x_1,y_1)(1_2,xy_2)(x_2,y_2)\\
 \gamma_1=(1_3,y_3)(x_3,xy_3)(1_4,y_4)(x_4,xy_4),
 & \ \   \gamma_4=(1_3,xy_3)(x_3,y_3)(1_4,xy_4)(x_4,y_4),\\
\gamma=(1_3,x_3)(xy_3,y_3)(1_4,x_4)(xy_4,y_4).
\end{align*}
Let $1\leq i,j\leq 4$. Since $\a\in Z(\Aut(\Gamma)_{(\G)})$, $R(G)=\langle \a,\beta_1\gamma_1\rangle=\langle\beta_4\gamma_4,\beta_1\gamma_1\rangle\cong\mz_2\times\mz_2$ is a semiregular normal subgroup of $\Aut(\Gamma)$ with orbit set $\G$, and $\langle \beta_1\gamma_1, \beta_4\gamma\rangle\cong\mz_2\times\mz_2$ is another semiregular subgroups of $\Aut(\Gamma)_{(\G)}$ with orbit set $\G$, which cannot be conjugate with $R(G)$ in $\Aut(\Gamma)$ because $R(G)\unlhd A$. Also, this can be easily checked by {\sc Magma}~\cite{magma}. By Theorem~\ref{PCI-BabaiSimilar}, $\Gamma$ is not $4$PCI, and hence $\mathbb{Z}_2\times\mathbb{Z}_2$ is not $4$PCI.

In what follows we construct a $4$-PCayley graph $\Gamma=\Cay(G,S_{i,j}: 1\leq i,j\leq 4)$ on each of the groups $\mz_4$, $\mz_6$ and $\mz_3\times \mz_3$. By {\sc Magma}~\cite{magma}, we may have $R(G)\unlhd \Aut(\Gamma)$, and there is another semiregular subgroup of $\Aut(\Gamma)$ isomorphic to $R(G)$ with orbit set $\G$, which cannot be conjugate to $R(G)$  in $\Aut(\Gamma)$. By Theorem~\ref{PCI-BabaiSimilar}, $\Gamma$ is not $4$PCI, and hence $\mz_4$, $\mz_6$ and $\mz_3\times \mz_3$, are not $4$PCI. These constructions are as follows: for
$G=\mathbb{Z}_4=\langle x\rangle$, take $S_{1,2}=S_{2,1}=S_{3,4}=S_{4,3}=\{x,x^3\}$,
$S_{1,3}=S_{2,4}=\{1,x\}$, $S_{3,1}=S_{4,2}=\{1,x^3\}$, $S_{2,3}=S_{3,2}=G$, and
$S_{i,j}=\emptyset$ for all other $1\leq i,j\leq 4$; for $G=\mathbb{Z}_6=\langle x\rangle$, take $S_{1,3}=\{x^2,x^4, x^5\}$, $S_{2,3}=S_{3,2}=\{1,x^2,x^4\}$, $S_{2,4}=\{1,x,x^4\}$, $S_{3,1}=\{x,x^2,x^4\}$, $S_{4,2}=\{1,x^2,x^5\}$, and $S_{i,j}=\emptyset$ for all other $1\leq m\leq 4$; for $G=\mathbb{Z}_3\times \mathbb{Z}_3=\langle x\rangle\times \langle y\rangle$, take $S_{1,3}=\{1,x,y,y^2,xy, x^2y\}$, $S_{1,4}=\{x,y^2, xy^2\}$, $S_{2,3}=\{y,x^2,y^2,xy,x^2y^2\}$, $S_{2,4}=\{x,y,xy^2,x^2y^2,x^2y\}$, $S_{3,1}=\{1,y,y^2,x^2, xy^2,x^2y^2\}$, $S_{3,2}=\{x,y,y^2,xy,$ $x^2y^2\}$, $S_{3,4}=\{x,y,y^2,x^2y^2\}$, $S_{4,1}=\{y,x^2,x^2y\}$, $S_{4,2}=\{y^2,x^2,xy^2,xy,yx^2\}$, $S_{4,3}=\{y,y^2,x^2,xy\}$ and $S_{i,j}=\emptyset$ for others $1\leq i,j\leq 4$.

Assume $G=\mz_3=\langle x\rangle$. Take
$S_{1,2}=S_{1,3}=S_{2,1}=S_{3,1}=S_{4,5}=S_{4,6}=S_{5,4}=S_{6,4}=\{1\}$,
$ S_{1,4}=S_{2,4}=S_{2,5}=S_{4,1}=S_{4,2}=S_{5,2}=\{1,x,x^2\}$,
$ S_{2,3}=S_{5,6}=\{x\}, S_{3,2}=S_{6,5}=\{x^2\}$, and
$ S_{i,j}=\emptyset$ for all other $1\leq i\leq j\leq 6$.
Let $\Gamma=\Cay(G,S_{i,j}: 1\leq i,j\leq 6)$. By {\sc Magma}~\cite{magma},  $\Aut(\Gamma)\cong (\mz_3 \times \mz_3)\rtimes \mz_{2}$,  $R(G)\unlhd \Aut(\Gamma)$, and there is another semiregular subgroup of $\Aut(\Gamma)$ isomorphic to $R(G)$ with orbit set $\G$, which cannot be conjugate to $R(G)$ in $\Aut(\Gamma)$. By Theorem~\ref{PCI-BabaiSimilar}, $\Gamma$ is not a $6$PCI-graph, and hence $\mathbb{Z}_3$ is not $6$PCI.

Assume $G=D_6$. By Theorem~\ref{subgroupsmPCI}, $D_6$ is not $6$PCI because $\mz_3$ is not $6$PCI.
\end{proof}

\medskip

By Lemma~\ref{non-PCIs},  $\mz_3$ and $D_6$ are not $6$PCI. However, it is much more difficult to prove that $\mz_3$ and $D_6$ are $m$PCI for every $m\leq 5$. This can be done for $m\leq 3$ by computation in {\sc Magma}~\cite{magma}, but not for $m=4$ and $5$ because it is beyond the computing power of our computer.
Our computation in {\sc Magma} proceeds as follows: Let $G=\mz_3$ or $D_6$ and let  $\Gamma=\Cay(G,S_{i,j}: 1\leq i,j\leq m)$ be an $m$-PCayley graph of $G$.
Each $S_{ij}$ can be any subset of $G$, offering $2^{\vert G\vert}$ choices.
Note that $S_{ji}=S_{ij}^{-1}$ as $\Gamma$ is a graph.
Then the set $\{S_{i,j}\ |\ 1\leq i,j\leq m\}$ has $ 2^{\frac{\vert G\vert m(m-1)}{2}}$ choices.
For each possible $S_{i,j}$, we construct the $m$-PCayley graph $\Gamma$ in {\sc Magma}  and compute its full automorphism groups $\mathrm{Aut}(\Gamma)$.
If $\mathrm{Aut}(\Gamma)$ has a semiregular subgroup which has the same orbit set as $R(G)$ and is isomorphic but not conjugate to $R(G)$, then $\Gamma$ is not $m$PCI by Theorem~\ref{PCI-BabaiSimilar}.

Notably, for $G = D_6$ and $m = 4$, there are $2^{36} \approx 6.8 \times 10^{10}$ choices for $\{ S_{i,j} : 1 \leq i,j \leq m \} $, making it challenging to verify that $D_6$ is $4$PCI using the above approach. Therefore, for $m=4$ and $m=5$, we seek a theoretical proof.

Let $\Gamma=\Cay(G,S_{i,j}: 1\leq i,j\leq m)$ be an $m$-PCayley graph of a finite group $G$. Recall that $\G=\{G_1,G_2,\cdots,G_m\}$. For convenience of statement, we need more notations for $\Gamma$.
\begin{itemize}
\item Denote by $\Gamma(S_{i,j})$ the $m$-PCayley graph deduced from $\Gamma$ by replacing $S_{i,j}$ by $G\backslash S_{i,j}$.
\item
An element in $\G$ is called {\em a block} of $\Gamma$, and $\Gamma$ is said to be {\em block connected} if for any two different blocks $G_i$ and $G_j$ in $\G$, there is a series of blocks $G_i=G_{m_1},G_{m_2},\cdots,G_{m_t}=G_j$ such that $[G_{m_k},G_{m_{k+1}}]_\Gamma$ has at least one edge for every $1\leq k\leq t-1$.
\item
A subgraph $\Sigma$ of $\Gamma$ is called {\em a block subgraph} if $V(\Sigma)$ is a union of some blocks of $\Gamma$ such that for any two distinct blocks $G_i$ and $G_j$ in $V(\Sigma)$, either $[G_i\cup G_j]_\Sigma=[G_i\cup G_j]_\Gamma$ or $[G_i\cup G_j]_\Sigma$ has no edge. Write $\G(\Sigma)=\{G_i\ |\ G_i\subseteq V(\Sigma)\}$. Then $|\G(\Sigma)|$ is called the {\em length} of $\Sigma$.

\item A block subgraph $\Sigma$ of $\Gamma$ is called {\em induced}, if $\Sigma$ is an induced subgraph of $\Gamma$, and a maximal induced block connected subgraph is called a {\em block component} of $\Gamma$. For a bock component $\Sigma$ of $\Gamma$, $\Aut(\Sigma)_{(\G)}$ can be viewed as a subgroup of  $\Aut(\Gamma)_{(\G)}$ by fixing $V(\Gamma)/V(\Sigma)$ pointwise.
\end{itemize}

From the above notation, it is easy to see the following proposition.

\begin{prop}\label{lm:blocksubgraphs}
Let $\Gamma=\Cay(G,S_{i,j}: 1\leq i,j\leq m)$ be an $m$-PCayley graph of a finite group $G$.
Let $\Aut(\Gamma)_{(\G)}$ be the subgroup of $\Aut(\Gamma)$ fixing $\G$ pointwise. Then
\begin{itemize}
\item[\rm (1)] $\Aut(\Gamma)_{(\G)}=\Aut(\Gamma(S_{i,j}))_{(\G)}$,   for any $1\leq i,j\leq m $ with $i\neq j$.
\item[\rm (2)] Let $\Sigma_1,\Sigma_2,\cdots,\Sigma_t$ be all block components of $\Gamma$. Then $\Aut(\Gamma)_{(\G)}=\Aut(\Sigma_1)_{(\G)}\times \Aut(\Sigma_2)_{(\G)}\times \cdots \times \Aut(\Sigma_t)_{(\G)}$.
\end{itemize}
\end{prop}

Transitive permutation groups $M$ of degree $6$ containing a regular dihedral subgroup will be repeatedly used later. If $M$ is primitive then it is $2$-transitive (see  \cite[Theorem 25.6]{HW}), and then $M=\S_6$ or $\S_5$ (see \cite{RoneyD}). If $M$ is imprimitive, then all non-trivial imprimitive blocks of $M$ have length $2$ or $3$, and hence $M\leq (\S_3\times\S_3)\rtimes\mz_2$ or $(\mz_2\times\mz_2\times\mz_2)\rtimes\S_3$. For the later, a Sylow $3$-subgroup of $M$ has order $3$, and so all subgroups of order $3$ are conjugate in $M$. Note that a semiregular subgroup of order $3$ of $\S_6$ uniquely determines a regular dihedral subgroup containing this subgroup of order $3$, and hence two regular dihedral subgroups of $M$ are conjugate if and only if their semiregular Sylow $3$-subgroups are conjugate in $M$. This implies that if a Sylow $3$-subgroup of $M$ has order $3$ then all regular dihedral subgroups of $M$ are conjugate. Assume that $3^2\mid |M|$. Then $M$ can be listed in the following Table~\ref{pgdegree6}, and this can be checked by {\sc Magma}~\cite{magma}.

\begin{table}[h]
\[
 \begin{array}{  |c | c|c|c|c|  } \hline
M & {\rm Primitive} & \mbox{$D_6$ Conjugate} & \mbox{ $D_6$ Normal} &  \mbox{ $D_6$ Unique}\\ \hline
 \S_6 & Yes & Yes & No & No \\ \hline
 (\S_3\times \S_3)\rtimes\mz_2 & No & Yes & No & No\\ \hline
 (\mz_3\times \mz_3)\rtimes\mz_2  & No & Yes & Yes & Yes\\ \hline
 (\mz_3\times \mz_3)\rtimes(\mz_2\times\mz_2) & No & No & Yes & No\\ \hline
 \end{array}
\]
\caption{Permutation groups $M$ of degree $6$ containing regular $D_6$ with $3^2\mid |M|$}\label{pgdegree6}
\end{table}

\begin{lem}\label{5-PCIs}
$\mz_3$ and $D_6$ are $m$PCI for every $m\leq 5$.
\end{lem}

\begin{proof}
Write $G=\mz_3$ or $D_6$. Let $\Gamma=\Cay(G,S_{i,j}: 1\leq i,j\leq m)$ be an $m$-PCayley graph of $G$ with $m\leq 5$, and $L\cong G$ a semiregular subgroup of $\Aut(\Gamma)$ with the same orbit set $\G=\{G_1,G_2,\cdots,G_m\}$ as $R(G)$.
Then $R(G),L\leq \Aut(\Gamma)_{(\G)}$, and by Theorem~\ref{PCI-BabaiSimilar}, it suffices to show that $R(G)$ and $L$ are conjugate in $\Aut(\Gamma)_{(\G)}$.

Since we only consider the group $\Aut(\Gamma)_{(\G)}$, by Proposition~\ref{lm:blocksubgraphs}~(1), we may assume $|S_{i,j}|\leq |G|/2$ for all $1\leq i,j\leq m$ with $i\not=j$. Let $\Sigma$ be a block component of $\Gamma$ with length $r$. Then $\Sigma$ is a $r$-PCayley graph of $G$  with $1\leq r\leq m\leq 5$.

\medskip
\noindent{\bf Case 1:} $G=\mz_3$.

Since $|S_{i,j}|\leq |G|/2=3/2$, the induced subgraph $[G_i,G_j]_\Gamma$ is either a matching or an empty graph, and since $\Sigma$ is a block component of $\Gamma$, $\Aut(\Sigma)_{(\G)}$ is faithful on each $G_i\in\G(\Sigma)$, forcing that $\Aut(\Sigma)_{(\G)}\cong \mz_3$ or $S_3$ as $|G_i|=3$. Note that $R(G)\cong L\cong \mz_3$. Then $R(G)^{V(\Sigma)}=L^{V(\Sigma)}$, the unique normal Sylow $3$-subgroup of  $\Aut(\Sigma)_{(\G)}$, where $R(G)^{V(\Sigma)}$ and $L^{V(\Sigma)}$ are the constituents of $R(G)$ and $L$ on $V(\Sigma)$, respectively, that is, the restricted permutation groups of $R(G)$ and $L$ on $V(\Sigma)$. In particular, if $\Gamma$ is block connected, then $R(G)=L$ and we are done. Thus, we assume that $\Gamma$ is not block connected, implying that $r\leq 4$ for $\Sigma$.

\medskip
\noindent{\bf Claim 1:} Let $\Delta$ be a union of some block components of $\Gamma$ except $\Sigma$. Assume that $R(G)^{V(\Delta)}$ and $L^{V(\Delta)}$ are conjugate in $\Aut(\Delta)_{(\G)}$ and that $\Aut(\Sigma)_{(\G)}\cong \S_3$. Then $R(G)^{V(\Delta\cup \Sigma)}$ and $L^{V(\Delta\cup \Sigma)}$ are conjugate in $\Aut(\Delta\cup \Sigma)_{(\G)}$.

Since $R(G)\cong L\cong\mz_3$, we may assume that
$R(G)^{V(\Delta\cup \Sigma)}=\langle \a_1\a_2\rangle$, where $\a_1\in \Aut(\Delta)_{(\G)}$ and  $\a_2\in \Aut(\Sigma)_{(\G)}$ have order $3$. By Proposition~\ref{lm:blocksubgraphs}~(2), we have $\Aut(\Delta)_{(\G)}\leq \Aut(\Delta\cup\Sigma)_{(\G)}$ and $\Aut(\Sigma)_{(\G)}\leq \Aut(\Delta\cup\Sigma)_{(\G)}$, where $\Aut(\Delta)_{(\G)}$ fixes $V(\Sigma)$ pointwise and $\Aut(\Sigma)_{(\G)}$ fixes $V(\Delta)$ pointwise. Since $R(G)^{V(\Delta)}$ and $L^{V(\Delta)}$ are conjugate in $\Aut(\Delta)_{(\G)}$, there is $\g\in \Aut(\Delta)_{(\G)}\leq \Aut(\Delta\cup\Sigma)_{(\G)}$ such that $(R(G)^{V(\Delta)})^\g=L^{V(\Delta)}$, and then $(\a_1\a_2)^\g=\a_1^\g\a_2$ and  $L^{V(\Delta)}=\langle \a_1^\g\rangle$. Since $\Aut(\Sigma)_{(\G)}$ has the unique Sylow $3$-subgroup $\langle\a_2\rangle$, we have $L^{V(\Sigma)}=\langle \a_2\rangle$. It follows $L^{V(\Delta\cup\Sigma)}=\langle \a_1^\g\a_2\rangle$ or $\langle \a_1^\g\a_2^2\rangle$.

If $L^{V(\Delta\cup\Sigma)}=\langle \a_1^\g\a_2\rangle$ then $(R(G)^{V(\Delta\cup\Sigma)})^\g=\langle \a_1^\g\a_2\rangle=L^{V(\Delta\cup\Sigma)}$, and we are done. If $L^{V(\Delta\cup\Sigma)}=\langle \a_1^\g\a_2^2\rangle$, then $\Aut(\Sigma)_{(\G)}\cong \S_3$ implies that there is $\b\in \Aut(\Sigma)_{(\G)}\leq \Aut(\Delta\cup\Sigma)_{(\G)}$ such that $\a_2^\b=\a_2^2$ and $(\a_1^\g)^\b=\a_1^\g$. Thus, $(R(G)^{V(\Delta\cup \Sigma)})^{\g\b}=\langle \a_1^\g\a_2\rangle^\b=\langle \a_1^\g\a_2^\b\rangle=\langle \a_1^\g\a_2^2\rangle=L^{V(\Delta\cup \Sigma)}$, as claimed.

\medskip
It is easy to see that if $r=1$ or $2$ then $\Aut(\Sigma)_{(\G)}\cong\S_3$. If every block component of $\Gamma$ has length $1$ or $2$, applying Claim~1 repeatedly we have that $R(G)$ and $L$ are conjugate in $\Aut(\Gamma)_{(\G)}$, as required. Thus, we may assume that $\Gamma$ has a block component of length at least $3$, say $\Sigma$. Recall that $R(G)^{V(\Sigma)}=L^{V(\Sigma)}$. Since $m\leq 5$, every  component of $\Gamma$ except $\Sigma$ has length $1$ or $2$, and Claim~1 implies that $R(G)$ and $L$ are conjugate in $\Aut(\Gamma)_{(\G)}$. This completes the proof of Case~1.

\medskip
\noindent{\bf Case 2:} $G=D_6$.

\medskip
Compared with Claim~1, the following is more general for $D_6$.

\medskip
\noindent{\bf Claim 2:} Let $\Delta$ be a union of some block components of $\Gamma$ except $\Sigma$. Assume $R(G)^{V(\Sigma)}$ and $L^{V(\Sigma)}$ are conjugate in $\Aut(\Sigma)_{(\G)}$, and $R(G)^{V(\Delta)}$ and $L^{V(\Delta)}$ are conjugate in $\Aut(\Delta)_{(\G)}$. Then $R(G)^{V(\Sigma\cup\Delta)}$ and $L^{V(\Sigma\cup\Delta)}$ are conjugate in $\Aut(\Sigma\cup\Delta)_{(\G)}$.

\medskip

Note that $R(G)\cong L\cong D_6$. We may assume  $R(G)^{V(\Sigma\cup\Delta)}=\langle \a_1\a_2,\b_1\b_2\rangle$, where $\a_1,\b_1\in \Aut(\Sigma)_{(\G)}\leq \Aut(\Sigma\cup\Delta)_{(\G)}$ and $\a_2,\b_2\in \Aut(\Delta)_{(\G)}\leq \Aut(\Sigma\cup\Delta)_{(\G)}$. Furthermore, $o(\a_1)=o(\a_2)=3$ and $o(\b_1)=o(\b_2)=2$ with $\a_1^{\b_1}=\a_1^{-1}$ and $\a_2^{\b_2}=\a_2^{-1}$.
Clearly,  $R(G)^{V(\Sigma)}=\langle \a_1,\b_1\rangle$ and $R(G)^{V(\Delta)}=\langle \a_2,\b_2\rangle$. By Proposition~\ref{lm:blocksubgraphs}~(2), $[\Aut(\Sigma)_{(\G)},\Aut(\Delta)_{(\G)}]=1$, that is, $\Aut(\Sigma)_{(\G)}$ commutes with $\Aut(\Delta)_{(\G)}$ pointwise.

By assumption, there are $\a\in \Aut(\Sigma)_{(\G)}\leq \Aut(\Sigma\cup\Delta)_{(\G)}$ and $\b\in \Aut(\Delta)_{(\G)}\leq \Aut(\Sigma\cup\Delta)_{(\G)}$ such that
 $(R(G)^{V(\Sigma)})^\a=L^{V(\Sigma)}$ and $(R(G)^{V(\Delta)})^{\b}=L^{V(\Delta)}$.
Note that $\a,\a_1,\b_1\in \Aut(\Sigma)_{(\G)}$ and $\b,\a_2,\b_2\in \Aut(\Delta)_{(\G)}$. Since $[\Aut(\Sigma)_{(\G)},\Aut(\Delta)_{(\G)}]=1$, we have $(R(G)^{V(\Sigma\cup\Delta)})^\a=\langle \a_1\a_2,\b_1\b_2\rangle^\a=\langle \a_1^\a\a_2,\b_1^\a\b_2\rangle$ and $(R(G)^{V(\Sigma\cup\Delta)})^\b=\langle \a_1\a_2,\b_1\b_2\rangle^\b=\langle \a_1\a_2^\b,\b_1\b_2^\b\rangle$. It follows
$$L^{V(\Sigma)}=\langle \a_1^\a,\b_1^\a\rangle\mbox{ and }
L^{V(\Delta)}=\langle \a_2^\b,\b_2^\b\rangle,$$ where $(\a_1^\a)^{\b_1^\a}=(\a_1^2)^\a$ and   $(\a_2^\b)^{\b_2^\b}=(\a_2^2)^\b$. Furthermore, $(R(G)^{V(\Sigma\cup\Delta)})^{\a\b}=\langle \a_1^\a\a_2^\b,\b_1^\a\b_2^\b\rangle$.

Note that $L^{V(\Sigma\cup\Delta)}\cong D_6$. Let $x\in L^{V(\Sigma\cup\Delta)}$ be an involution. Then $x=ab$ with $a\in \{\b_1^\a,\b_1^\a\a_1^\a,\b_1^\a(\a_1^2)^\a\}$ and $b\in  \{\b_2^\b,\b_2^\b\a_2^\b,\b_2^\b(\a_2^2)^\b\}$, and $L^{V(\Sigma\cup\Delta)}$ has a unique Sylow $3$-subgroup $L_3$ with $L_3=\langle \a_1^\a\a_2^\b\rangle$ or $\langle \a_1^\a(\a_2^2)^\b\rangle$. It follows that there is $y\in L_3$ such that $x^y=\b_1^\a b^y\in L^{V(\Sigma\cup\Delta)}$, where $b^y\in  \{\b_2^\b,\b_2^\b\a_2^\b,\b_2^\b(\a_2^2)^\b\}$, and $L^{V(\Sigma\cup\Delta)}=\langle L_3,x^y\rangle$. For $L_3=\langle \a_1^\a\a_2^\b\rangle$, we have $L^{V(\Sigma\cup\Delta)}=\langle \a_1^\a\a_2^\b, \b_1^\a\b_2^\b\rangle, \langle \a_1^\a\a_2^\b, \b_1^\a\b_2^\b\a_2^\b\rangle$ or $\langle \a_1^\a\a_2^\b, \b_1^\a\b_2^\b(\a_2^2)^\b\rangle$, and by taking $\g=1, (\a_2^2)^\b$ or $\a_2^\b$ respectively, we have $(R(G)^{V(\Sigma\cup\Delta)})^{\a\b\g}=\langle \a_1^\a\a_2^\b,\b_1^\a\b_2^\b\rangle^\g=L^{V(\Sigma\cup\Delta)}$. For $L_3=\langle \a_1^\a(\a_2^2)^\b\rangle$, we have $L^{V(\Sigma\cup\Delta)}=\langle \a_1^\a(\a_2^2)^\b, \b_1^\a\b_2^\b\rangle, \langle \a_1^\a(\a_2^2)^\b, \b_1^\a\b_2^\b\a_2^\b\rangle$ or $\langle \a_1^\a(\a_2^2)^\b, \b_1^\a\b_2^\b(\a_2^2)^\b\rangle$, and by taking $\g=\b_2^\b, \b_2^\b(\a_2^2)^\b$ or $\b_2^\b\a_2^\b$ respectively, we have $(R(G)^{V(\Sigma\cup\Delta)})^{\a\b\g}=\langle \a_1^\a\a_2^\b,\b_1^\a\b_2^\b\rangle^\g=L^{V(\Sigma\cup\Delta)}$, as claimed.

\medskip

An edge in an induced subgraph $[G_i\cup G_j]_\Gamma$ is called a {\em $K_{3,3}$-edge} if $[G_i\cup G_j]_\Gamma\cong 2K_{3,3}$, and a {\em $K_{2,2}$-edge} if $[G_i\cup G_j]_\Gamma\cong 3K_{2,2}$. Write $G=\langle a,b\ |\ a^3=b^2=1,bab=a^2\rangle$ and $H=\langle a\rangle$. Then $G_i=H_i\cup (bH)_i$ for $1\leq i\leq m$, where $H_i=\{1_i,a_i,a^2_i\}$ and $(bH)_i=\{b_i,(ba)_i,(ba^2)_i\}$.

Recall that $|S_{i,j}|\leq |G|/2=3$. It is easy to see that if $|S_{i,j}|=1$ then $[G_i\cup G_j]_\Gamma\cong 6K_2$; if  $|S_{i,j}|=2$ then $[G_i\cup G_j]_\Gamma\cong 3K_{2,2}$ or $2(K_{3,3}-3K_2)$ (two $6$-cycles); if $|S_{i,j}|=3$ then $[G_i\cup G_j]_\Gamma\cong 2K_{3,3}$ or $\Aut([G_i\cup G_j])_{(\G)}$ is faithful on each of $G_i$ and $G_j$. These facts can be obtained from {\sc Magma}~\cite{magma} or the fact that there is no connected symmetric or semisymmetric cubic graphs of order $12$ (see \cite{CD,CMMP}). Noting that $R(H)$ induces a group of automorphisms of $[G_i\cup G_j]$ fixing $H_i$ and $H_j$, we have the following simple observation.

\medskip
\noindent{\bf Observation:} Assume that $[G_i\cup G_j]\not\cong 3K_{2,2}$, $12K_1$. If $\a\in \Aut([G_i\cup G_j])_{(\G)}$ fixes $H_i$ setwise then $\a$ fixes $H_j$ setwise. If further $[G_i\cup G_j]\not\cong 2K_{3,3}$, then the restriction of $\Aut([G_i\cup G_j])_{(\G)}$ on each of $G_i$ and $G_j$ is faithful.

\medskip

\medskip Let $\Lambda$ be a block subgraph of $\Gamma$. Then  $R(G)^{V(\Lambda)},L^{V(\Lambda)}\leq \Aut(\Lambda)_{(\G)}$.
Write $$L(H)=\{ H_i\ |\ G_i\in\G(\Lambda),\{H_i,(bH)_i\} \mbox{ is an imprimitive block system of } L \mbox{ on } {G_i}\}.$$
Denote by $\Aut(\Lambda)_{(\G)}^+$ the subgroup of $\Aut(\Lambda)_{(\G)}$ fixing $H_i$ setwsie for every $H_i\in L(H)$ and denote by $R_3$ and $L_3$ the unique Sylow $3$-subgroups of $R(G)^{V(\Lambda)}$ and $L^{V(\Lambda)}$, respectively. Clearly, $\Aut(\Lambda)_{(\G)}^+$ are dependent on the choice of $L$ and $L_3\leq  \Aut(\Lambda)_{(\G)}^+$.
Since $R(G)$ has the block system $\{H_i,(bH)_i\}$ for every $G_i\in\G(\Lambda)$, we have $R_3\leq \Aut(\Lambda)_{(\G)}^+$.

Denote by $\Lambda^3$ the block subgraph of $\Lambda$ by deleting all $K_{3,3}$-edges from $\Lambda$. Then $\Lambda^3$ is also a block subgraph of $\Gamma$.
If $[G_i,G_j]_\Gamma\cong 2K_{3,3}$ with $G_i,G_j\in\G(\Lambda)$, then $S_{i,j}=H$ or $bH$, and then $H_i$ and $H_j$ are blocks of $L$ on $G_i$ and $G_j$, respectively. Thus,
\begin{equation}\label{Lambda^3Auto}
\Aut(\Lambda)_{(\G)}^+=\Aut(\Lambda^3)_{(\G)}^+.
\end{equation}

\medskip
\noindent{\bf Claim 3:} Assume that $\Lambda^3$ is block connected. Then $R(G)^{V(\Lambda)}$ and $L^{V(\Lambda)}$ are conjugate in $\Aut(\Lambda)_{(\G)}$, and further in  $\Aut(\Lambda)_{(\G)}^+$, unless $\Aut(\Lambda)_{(\G)}$ is isomorphic to $\S_6$ and is faithful on every $G_i\in\G(\Lambda)$.

\medskip
First assume that $\Lambda^3$ contains no $K_{2,2}$-edges. Since $\Lambda^3$ is block connected, Observation implies that $\Aut(\Lambda)_{(\G)}$ is faithful on every $G_i\in\G(\Lambda)$, and hence $\Aut(\Lambda)_{(\G)}\lesssim\S_6$. Let $S$ be a Sylow $3$-subgroup of $\Aut(\Lambda)_{(\G)}^+$. Since $R_3,L_3\leq \Aut(\Lambda)_{(\G)}^+$, we have $|S|=3$ or $9$.

Let $|S|=3$. By Sylow theorem,  there is $\a\in \Aut(\Lambda)_{(\G)}^+$ such that $R_3^\a=L_3$.
Note that $|G_i|=6$, and hence for a given semiregular subgroup of order $3$ in the symmetric group $\S_{G_i}$, there is a unique regular dihedral subgroup of $\S_{G_i}$ containing the semiregular subgroup. Thus, $(R(G)^{G_i})^\a=L^{G_i}$ as $R_3^\a=L_3$, and hence  $(R(G)^{V(\Lambda)})^\a=L^{V(\Lambda)}$, because $\Aut(\Lambda)_{(\G)}$ has a faithful restriction on $G_i$.

Let $|S|=9$. Since $\Aut(\Lambda)_{(\G)}$ is faithful on $G_i$, by Table~\ref{pgdegree6} we have $\Aut(\Lambda)_{(\G)}\cong(\Aut(\Lambda)_{(\G)})^{G_i} \cong\S_6,(\S_3\times\S_3)\rtimes\mz_2,(\mz_3\times\mz_3)\rtimes\mz_2$ or $(\mz_3\times\mz_3)\rtimes(\mz_2\times\mz_2)$.

Suppose $\Aut(\Lambda)_{(\G)}\cong\S_6$. By Theorem~\ref{SRGconj}, $R(G)^{G_i}$ and $L^{G_i}$ are conjugate in $(\Aut(\Lambda)_{(\G)})^{G_i}$, and hence $R(G)^{V(\Lambda)}$ and $L^{V(\Lambda)}$ are conjugate in $\Aut(\Lambda)_{(\G)}$, but may not in $\Aut(\Lambda)_{(\G)}^+$ . In this case, $\Aut(\Lambda)_{(\G)}$ is faithful on every $G_i\in \G(\Lambda)$. Suppose  $\Aut(\Lambda)_{(\G)}\cong(\S_3\times\S_3)\rtimes\mz_2$. Then $\Aut(\Lambda)_{(\G)}^+\cong \S_3\times\S_3$ and so all semiregular subgroups of order $3$ of $\Aut(\Lambda)_{(\G)}^+$ are conjugate because all semiregular subgroups of order $3$ of $(\Aut(\Lambda)_{(\G)}^+)^{G_i}$ are conjugate. It follows that $R(G)^{V(\Lambda)}$ and $L^{V(\Lambda)}$ are conjugate in $\Aut(\Lambda)_{(\G)}^+$. Suppose $\Aut(\Lambda)_{(\G)}\cong (\mz_3\times\mz_3)\rtimes\mz_2$. Then $\Aut(\Lambda)_{(\G)}^+\cong \mz_3\times\mz_3$ and $(\Aut(\Lambda)_{(\G)})^{G_i}$ has a unique regular dihedrant on $G_i$. Thus, $\Aut(\Lambda)_{(\G)}$  has a unique semiregular dihedrant, and hence $R(G)^{V(\Lambda)}=L^{V(\Lambda)}$. At last suppose  $\Aut(\Lambda)_{(\G)}\cong(\mz_3\times\mz_3)\rtimes(\mz_2\times\mz_2)$. Then $\Aut(\Lambda)_{(\G)}^+\cong(\mz_3\times\mz_3)\rtimes\mz_2$ and $(\Aut(\Lambda)_{(\G)})^{G_i}$ has exactly two distinct normal regular subgroups isomorphic to $D_6$ on $G_i$. Suppose $R(G)^{V(\Lambda)}\not=L^{V(\Lambda)}$. Then $R(G)^{V(\Lambda)}$ and $L^{V(\Lambda)}$ are the two distinct normal semiregular subgroups isomorphic to $D_6$ in $\Aut(\Lambda)_{(\G)}$. Thus, $\Aut(\Lambda)_{(\G)}=R(G)^{V(\Lambda)}\times L^{V(\Lambda)}$. Let $\langle\a\rangle$ and $\langle\b \rangle$ be the unique Sylow $3$-subgroups of $R(G)^{V(\Lambda)}$ and $L^{V(\Lambda)}$, respectively. Then $\langle\a\rangle\times\langle \b\rangle$ is the unique Sylow $3$-subgroup of $\Aut(\Lambda)_{(\G)}$. Clearly, $\langle\a\rangle\leq R(G)^{V(\Lambda)}$ has orbit set $\{H_i,(bH)_i\}$ on $G_i$, which is also the orbit set of $\langle\b\rangle$. Since $\Aut(\Lambda)_{(\G)}$ is faithful on $G_i$, $\langle\a\b\rangle$ is transitive on one of $H_i$ and $(bH)_i$, say $H_i$, and hence fixes $(bH)_i$ pointwise. This forces that $\langle\a\b^2\rangle$ is transitive on $(bH)_i$, and fixes $H_i$ pointwise. By the arbitrary of $G_i\in \G(\Lambda)$, we have $S_{i,j}=\emptyset,H$ or $bH$ for any $G_i,G_j\in\G(\Lambda)$ with $i\not=j$. It is easy to check that $\g: x_i\mapsto (x^{-1})_i$ for all $x_i\in G_i\in \G(\Lambda)$ is an automorphism of $\Lambda$, and hence $\g\in\Aut(\Lambda)_{(\G)}$. Moreover, $\g$ maps $\a$ to $\b$ or $\b^2$, contradicting that $\Aut(\Lambda)_{(\G)}=R(G)^{V(\Lambda)}\times L^{V(\Lambda)}$. It follows that $R(G)^{V(\Lambda)}=L^{V(\Lambda)}$. This yields that Claim~3 is true when $\Lambda^3$ contains no $K_{2,2}$-edges.

\medskip
Now assume that $\Lambda^3$ contains $K_{2,2}$-edges. Let $G_i,G_j\in \G(\Lambda)$ such that  $[G_i,G_j]_{\Lambda^3}=T_1\cup T_2\cup T_3\cong 3K_{2,2}$ with $T_i\cong K_{2,2}$. Then $[G_i,G_j]_{\Lambda}=T_1\cup T_2\cup T_3$. Considering the kernel $K$ of $\Aut(\Lambda)_{(\G)}$ on $\{T_1,T_2,T_3\}$, we have $\Aut(\Lambda)_{(\G)}/K\lesssim \S_3$. Every orbit of $K$ on $G_i\cup G_j$ has length $1$ or $2$, and hence $K/K_1$ is a $2$-group, where $K_1$ is the kernel of $K$ on $G_i\cup G_j$, that is, the subgroup of $K$ fixing every vertex in $G_i\cup G_j$. Further let $G_k\in \G(\Lambda)$ with $k\not=i,j$ such that $[G_i,G_k]_{\Lambda^3}$ has at least one edge. If $[G_i,G_k]_{\Lambda^3}\not\cong 3K_{2,2}$, then Observation implies that $K_1$ fixes $G_k$ pointwise. If $[G_i,G_k]_{\Lambda^3}\cong 3K_{2,2}$, then the above argument implies that $K_1/K_2$ is a $2$-group, where $K_2$ is the kernel of $K_1$ on $G_k$. This implies that $K/K_3$ is a $2$-group, where $K_3$ fixes $G_s$ pointwise for all $G_s\in \G(\Lambda)$ such that $[G_i,G_s]_{\Lambda^3}$ has at leat one edge. Since $\Lambda^3$ is block connected, an induction on the distance between $G_i$ and $G_\ell$ in $\Lambda^3$ with $G_\ell\in \G(\Lambda)$ implies that $K$ is a $2$-group. It follows that a Sylow $3$-subgroup of $\Aut(\Lambda)_{(\G)}$ has order $3$. Since $R_3,L_3\in \Aut(\Lambda)_{(\G)}^+$, by the Sylow theorem there is $\a\in \Aut(\Lambda)_{(\G)}^+$ such that $R_3^\a=L_3$.

Recall that $[G_i,G_j]_\Lambda=T_1\cup T_2\cup
T_3\cong 3K_{2,2}$. For not making the notation too cumbersome, we  write $G_i=\{1,2,3,4,5,6\}$ and $G_j=\{1',2',3',4',5',6'\}$ with $V(T_1)=\{1,4,1',4'\}$, $V(T_2)=\{2,5,2',5'\}$ and $V(T_3)=\{3,6,3',6'\}$. Then  $[G_i,G_j]_\Lambda$ can be drawn as the following Figure~\ref{figure2}.

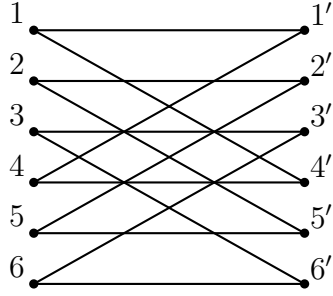
\begin{figure}[htp!]
  \centering
\begin{tikzpicture}[scale=0.45]
\draw [thick] (7,4)--(15,-0.5)--(7,-0.5)--(15,4)--(7,4);
\draw [thick] (7,2.5)--(15,2.5)--(7,-2)--(15,-2)--(7,2.5);
\draw [thick] (7,1)--(15,1)--(7,-3.5)--(15,-3.5)--(7,1);
\filldraw[thick,fill=black] (7,4) circle (3pt);
\node at (6.5,4.5) {$1$};
\filldraw[thick,fill=black] (7,2.5) circle (3pt);
\node at (6.5,3) {$2$};
\filldraw[thick,fill=black] (7,1) circle (3pt);
\node at (6.5,1.5) {$3$};
\filldraw[thick,fill=black] (7,-0.5) circle (3pt);
\node at (6.5,0) {$4$};
\filldraw[thick,fill=black] (7,-2) circle (3pt);
\node at (6.5,-1.5) {$5$};
\filldraw[thick,fill=black] (7,-3.5) circle (3pt);
\node at (6.5,-3) {$6$};
\filldraw[thick,fill=black] (15,4) circle (3pt);
\node at (15.5,4.5) {$1'$};
\filldraw[thick,fill=black] (15,2.5) circle (3pt);
\node at (15.5,3) {$2'$};
\filldraw[thick,fill=black] (15,1) circle (3pt);
\node at (15.5,1.5) {$3'$};
\filldraw[thick,fill=black] (15,-0.5) circle (3pt);
\node at (15.5,0) {$4'$};
\filldraw[thick,fill=black] (15,-2) circle (3pt);
\node at (15.5,-1.5) {$5'$};
\filldraw[thick,fill=black] (15,-3.5) circle (3pt);
\node at (15.5,-3) {$6'$};
\end{tikzpicture}
\caption{$[G_i,G_j]_\Lambda$ of $\Lambda^3$ contains $K_{2,2}$-edges}
\label{figure2}
\end{figure}
Let $R_3^\a=L_3=\langle\d\rangle$. Then $\d$ permutes the three $K_{2,2}$-copies of $[G_i,G_j]_\Lambda$ cyclicly, and without loss of any generality, we may assume that $\d^{G_i\cup G_j}=(1\ 2\ 3)(4\ 5\ 6)(1'\ 2'\ 3')(4'\ 5'\ 6')$. There is a unique dihedral subgroup of $\Aut([G_i,G_j]_\Lambda)_{(\G)}$  that is regular on $G_i$ and $G_j$, which in fact is $\langle \d, (1\ 4)(2\ 6)(3\ 5)(1'\ 4')(2'\ 6')(3'\ 5')\rangle$. This implies that $(R(G)^{G_i\cup G_j})^\a=L^{G_i\cup G_j}$. Write $\Re=\cup_{[G_i,G_j]_\Lambda\cong 3K_{2,2} \mbox{ with }G_i,G_j\in\G(\Lambda)} G_i$. Then $(R(G)^\Re)^\a=L^\Re$. Since $\Lambda^3$ is block connected, $\Aut(\Lambda)_{(\G)}$ is faithful on $\Re$ by Observation, and then $(R(G)^{V(\Lambda)})^\a=L^{V(\Lambda)}$. This completes the proof of Claim~3.

\medskip
Recall that $\Sigma$ is a block component of $\Gamma$ and has length $r$ with $1\leq r\leq m\leq 5$. Assume that $\Sigma^3$ is block connected. By Claim~3, $R(G)^{V(\Sigma)}$ and $L^{V(\Sigma)}$ are conjugate in $\Aut(\Sigma)_{(\G)}$.

Now assume that $\Sigma^3$ is not block connected. Let $\Sigma_1,\Sigma_2,\cdots,\Sigma_t$ be all block components of $\Sigma^3$. Then $t\geq 2$. By Eq~\ref{Lambda^3Auto}, $\Aut(\Sigma)_{(\G)}^+=\Aut(\Sigma^3)_{(\G)}^+=\Aut(\Sigma_1\cup\cdots\cup\Sigma_t)_{(\G)}^+$.
Clearly, $\Aut(\Sigma_i)_{(\G)}^+$ can be extended to a subgroup of $\Aut(\Sigma)_{(\G)}^+$ by fixing $V(\Sigma_j)$ pointwise for every $j\not=i$, and then we have $\Aut(\Sigma_i)_{(\G)}^+\leq \Aut(\Sigma)_{(\G)}^+$. Similar to Proposition~\ref{lm:blocksubgraphs}~(2), we have

\begin{equation}\label{component+product}
\Aut(\Sigma)_{(\G)}^+=\Aut(\Sigma_1)_{(\G)}^+\times \Aut(\Sigma_2)_{(\G)}^+\times \cdots \times \Aut(\Sigma_t)_{(\G)}^+.
\end{equation}

\medskip
 Let $X$ be a block component of $\Sigma^3$. Then $X$ is the subgraph of $[V(X)]_\Gamma$ by deleting all $K_{3,3}$-edges, and hence  $\Aut([V(X)]_\Gamma)_{(\G)}^+=\Aut(X)_{(\G)}^+$ by Eq~(\ref{Lambda^3Auto}).
 Since $t\geq 2$, there are $G_i\in \G(X)$ and $G_j\not\in \G(X)$ with $G_j\in\G(\Sigma)$ such that $[G_i,G_j]_\Gamma\cong 2K_{3,3}$. Then $\{H_i,(bH)_i\}$ is an imprimitive block system of $L$ on $G_i$ and hence $\Aut(X)_{(\G)}$ cannot be primitive on $G_i$. In particular, $\Aut(X)_{(\G)}$ cannot be faithful on $G_i$ and isomorphic to $\S_6$. By Claim~3, $R(G)^{V(X)}$ and $L^{V(X)}$ are conjugate in $\Aut(X)_{(\G)}^+$. Similar to Claim~2, we have the following result.

\medskip
\noindent{\bf Claim 4:} Let $Y$ be a union of some block components of $\Sigma^3$ except $X$. Assume that $R(G)^{V(Y)}$ and $L^{V(Y)}$ are conjugate in $\Aut(Y)_{(\G)}^+$, and
$\Aut(X)_{(\G)}^+$ has an involution $\delta$ normalizing $L^{V(X)}$ and inversing its elements of order $3$. Then $R(G)^{V(X\cup Y)}$ and $L^{V(X\cup Y)}$ are conjugate in $\Aut(X\cup Y)_{(\G)}^+$.

Since $R(G)\cong L\cong D_6$, we may assume  $R(G)^{V(X\cup Y)}=\langle \a_1\a_2,\b_1\b_2\rangle$, where $\a_1\in \Aut(X)_{(\G)}^+$, $\b_1\in \Aut(X)_{(\G)}$, $\a_2\in \Aut(Y)_{(\G)}^+$, $\b_2\in \Aut(Y)_{(\G)}$.  Then $R(G)^{V(X)}=\langle \a_1,\b_1\rangle$ and  $R(G)^{V(Y)}=\langle \a_2,\b_2\rangle$. Furthermore, $o(\a_1)=o(\a_2)=3$ and $o(\b_1)=o(\b_2)=2$ with $\a_1^{\b_1}=\a_1^{-1}$ and $\a_2^{\b_2}=\a_2^{-1}$. Since $R(G)^{V(X)}$ and $L^{V(X)}$ are conjugate in $\Aut(X)_{(\G)}^+$, there is $\a\in \Aut(X)_{(\G)}^+$ such that $(R(G)^{V(X)})^\a=L^{V(X)}$. Since $R(G)^{V(Y)}$ and $L^{V(Y)}$ are conjugate in $\Aut(Y)_{(\G)}^+$ by assumption, there is $\b\in \Aut( Y)_{(\G)}^+$ such that $(R(G)^{V( Y)})^{\b}=L^{V( Y)}$. Thus,
 $$L^{V(X)}=\langle \a_1^\a,\b_1^\a\rangle\mbox{ and }
L^{V( Y)}=\langle \a_2^\b,\b_2^\b\rangle,$$
where $\b_1^\a\in\Aut(X)_{(\G)}$, $\b_2^\b\in\Aut(Y)_{(\G)}$, $\a_1^\a\in\Aut(X)_{(\G)}^+$ and $\a_2^\b\in\Aut(Y)_{(\G)}^+$. Since $L^{V(X\cup Y)}\cong D_6$, $L^{V(X\cup Y)}$ has a unique Sylow $3$-subgroup, say $L_3$. Then $L_3=\langle \a_1^\a\a_2^\b\rangle$ or $\langle \a_1^\a(\a_2^2)^\b\rangle$.

Since $\a$ fixes $V(Y)$ pointwise and $\b$ fixes $V(X)$ pointwise, we have
$$(R(G)^{V(X\cup Y)})^{\a\b}=\langle \a_1\a_2,\b_1\b_2\rangle^{\a\b}=\langle\a_1^\a\a_2,\b_1^\a\b_2\rangle^\b=
\langle\a_1^\a\a_2^\b,\b_1^\a\b_2^\b\rangle.$$

Let $x\in L^{V(X\cup Y)}$ be an involution. Then $x=ab$ with $a\in \{\b_1^\a,\b_1^\a\a_1^\a,\b_1^\a(\a_1^2)^\a\}$ and $b\in  \{\b_2^\b,\b_2^\b\a_2^\b,\b_2^\b(\a_2^2)^\b\}$. Since $L^{V(X\cup Y)}$ has the unique Sylow $3$-subgroup $L_3$ with $L_3=\langle \a_1^\a\a_2^\b\rangle$ or $\langle \a_1^\a(\a_2^2)^\b\rangle$, there is $y\in L_3$ such that $x^y=\b_1^\a b^y\in L^{V(X\cup Y)}$ with $b^y\in  \{\b_2^\b,\b_2^\b\a_2^\b,\b_2^\b(\a_2^2)^\b\}$, and $L^{V(X\cup Y)}=\langle L_3,x^y\rangle$. For $L_3=\langle \a_1^\a\a_2^\b\rangle$, we have $L^{V(X\cup Y)}=\langle \a_1^\a\a_2^\b, \b_1^\a\b_2^\b\rangle, \langle \a_1^\a\a_2^\b, \b_1^\a\b_2^\b\a_2^\b\rangle$ or $\langle \a_1^\a\a_2^\b, \b_1^\a\b_2^\b(\a_2^2)^\b\rangle$, and by taking $\g=1, (\a_2^2)^\b$ or $\a_2^\b$ respectively, we have $(R(G)^{V(X\cup Y)})^{\a\b\g}=\langle \a_1^\a\a_2^\b,\b_1^\a\b_2^\b\rangle^\g=L^{V(X\cup Y)}$. By Eq~\ref{component+product}, $\a\b\g\in\Aut(X\cup Y)_{(\G)}^+$. For $L_3=\langle \a_1^\a(\a_2^2)^\b\rangle$, we have $L^{V(X\cup Y)}=\langle \a_1^\a(\a_2^2)^\b, \b_1^\a\b_2^\b\rangle, \langle \a_1^\a(\a_2^2)^\b, \b_1^\a\b_2^\b\a_2^\b\rangle$ or $\langle \a_1^\a(\a_2^2)^\b, \b_1^\a\b_2^\b(\a_2^2)^\b\rangle$. Recall that $L^{V(X)}=\langle \a_1^\a,\b_1^\a\rangle$. By assumption, $(\a_1^\a)^\d=(\a_1^\a)^{-1}$ and $(L^{V(X)})^\d=L^{V(X)}$, where $\d\in \Aut(X)_{(\G)}^+$ is an involution. Since $L^{V(X)}$ has three involutions, $\d$ fixes one of them, and we may let  $(\b_1^\a)^\d=\b_1^\a$. Since $\d$ fixes $V(Y)$ pointwise,  we have $(\b_1^\a\b_2^\b)^\d=\b_1^\a\b_2^\b$ and $(\a_2^\b)^\d=\a_2^\b$. Then $(R(G)^{V(X\cup Y)})^{\a\b\d}=\langle \a_1^\a\a_2^\b,\b_1^\a\b_2^\b\rangle^{\d}=\langle (\a_1^\a)^2\a_2^\b,\b_1^\a\b_2^\b\rangle=\langle \a_1^\a(\a_2^2)^\b,\b_1^\a\b_2^\b\rangle$.
By taking $\g=1, (\a_2^2)^\b$ or $\a_2^\b$ respectively, we have $(R(G)^{V(X\cup Y)})^{\a\b\d\g}=L^{V(X\cup Y)}$, where $\a\b\d\g\in \Aut(X\cup Y)_{(\G)}^+$. This completes the proof of Claim~4.

\medskip
Now we are ready to finish the proof. If $R(G)^{V(\Lambda)}$ and $L^{V(\Lambda)}$ are conjugate in $\Aut(\Lambda)_{(\G)}$ for all components $\Lambda$ of $\Gamma$, then Claim~2 implies that $R(G)$ and $L$ are conjugate in $\Aut(\Gamma)$ and we are done. Recall that $\Sigma$ is an arbitrary block component of $\Gamma$ with length $1\leq r\leq m\leq 5$. To finish the proof, by Claim~2 we only need to show that $R(G)^{V(\Sigma)}$ and $L^{V(\Sigma)}$ are conjugate in $\Aut(\Sigma)_{(\G)}$.

If $\Sigma^3$ is block connected, then Claim~3 implies that $R(G)^{V(\Sigma)}$ and $L^{V(\Sigma)}$ are conjugate in $\Aut(\Sigma)_{(\G)}$, and we are done. Thus, we may assume that $\Sigma^3$ is not block connected. Let $X$ be a block component of $\Sigma^3$. Then there exist $G_i\in\G(X)$ and $G_j\in\G(\Sigma)$ with $G_j\not\in\G(X)$ such that $[G_i,G_j]_\Gamma\cong 2K_{3,3}$. Thus, $\{H_i,(bH)_i\}$ is an imprimitive block system of $L$ on $G_i$. By Claim~3, $R(G)^{V(X)}$ and $L^{V(X)}$ are conjugate in $\Aut(X)_{(\G)}^+$.

Suppose $X$ has length $1$. Then $\G(X)=\{G_i\}$, $V(X)=G_i$ and $\Aut(X)_{(\G)}^+\cong\S_3\times\S_3$. It is easy to see that $\Aut(X)_{(\G)}^+$ has an involution fixing $L^{V(X)}$ and inversing its elements of order $3$. Suppose $X$ has length $2$. We may write $\G(X)=\{G_i,G_k\}$. If $X\not\cong 3K_{2,2}$, by Observation we have that $\Aut(X)_{(\G)}^+$ is faithful on $G_i$.  Then it is easy to see that $\Aut(X)_{(\G)}^+\cong\S_3\times\S_3$, and $\Aut(X)_{(\G)}^+$ has an involution fixing $L^{V(X)}$ and inversing its elements of order $3$. If $X\cong 3K_{2,2}$, by {\sc Magma}~\cite{magma} we have $\Aut(X)^{+}_{\G}\cong \S_3$, and then it is easy to see that $\Aut(X)_{(\G)}^+$ has an involution fixing $L^{V(X)}$ and inversing its elements of order $3$ (note that $H_k$ may not be a block of $L$ on $G_k$).

If all block components of $\Sigma^3$ has length $1$ or $2$, then Claim~4, together with Claim~3 and the above paragraph, implies that $R(G)^{V(\Sigma)}$ and $L^{V(\Sigma)}$ are conjugate in $\Aut(\Sigma)_{(\G)}$, and we are done. Thus, we may assume that $X$ has length at least $3$. Recall that $R(G)^{V(X)}$ and $L^{V(X)}$ are conjugate in $\Aut(X)_{(\G)}^+$. Since $r\leq 5$, all other block components of $\Sigma^3$ except $X$ have length $1$ or $2$. Again by Claim~4,  $R(G)^{V(\Sigma)}$ and $L^{V(\Sigma)}$ are conjugate in $\Aut(\Sigma)_{(\G)}$, as required.
\end{proof}

Now we are ready to classify $m$PCI- and $m$PDCI-groups for $m\geq 4$.

\begin{theorem}\label{mPCI-groups} Let $m\geq 4$ be a positive integer and let $G$ be a finite group. Then
\begin{enumerate}
    \item[\rm (1)] $G$ is $m$PCI if and only if $G=\mz_1$, $\mathbb{Z}_2$, or $G=\mathbb{Z}_3$ or $D_6$ with $m=4$ or $5$;
    \item[\rm (2)] $G$ is $m$PDCI if and only if $G=\mz_1$ or $\mathbb{Z}_2$.
\end{enumerate}
\end{theorem}
\begin{proof}
 Let $\Gamma$ be an $m$-PCayley (di)graph of $G$,  If $G=\mz_1$ or $\mz_2$, then  $\Aut(\Gamma)$ has a unique semiregular subgroup with orbit set $\G=\{G_1,G_2,\cdots,G_m\}$, and hence $G$ is an $m$PCI-group and an $m$PDCI-group for any positive integer $m$. To prove parts~(1) and (2), we may assume that $|G|\geq 3$.

The sufficiency of part~(1) follows from Lemma~\ref{5-PCIs}, and to prove the necessity, let $G$ be an $m$PCI-group. By Lemma~\ref{pci2346}, $G$ is a $\{2,3\}$-group and so solvable. Let $N$ be a minimal normal subgroup of $G$. Then $N\cong\mz_2^r$ or $\mz_3^t$ for some positive integers $r$ and $t$. By Theorem~\ref{subgroupsmPCI}, $N$ is $m$PCI, and by Theorem~\ref{mtom-1}, any subgroup of $G$ is $4$PCI and so $N$ is $4$PCI.

Suppose $N\cong\mz_2^r$. If $r\geq 2$, then $\mz_2^2$ is $4$PCI, contradicting Lemma~\ref{non-PCIs}. Thus, $N=\mz_2$ and hence $N\leq Z(G)$, the center of $G$. If $G$ is a $2$-group, then $G$ has a subgroup isomorphic to $\mz_2^2$ or $\mz_4$ as $|G|\geq 3$, and hence $\mz_2^2$ or $\mz_4$ is $4$PCI, contradicting Lemma~\ref{non-PCIs}. Thus, $3\mid |G|$ and so $G$ has a subgroup isomorphic to $\mz_6$, implying $\mz_6$ is $4$PCI,  contradicting Lemma~\ref{non-PCIs}.

Thus, $N\cong\mz_3^t$. If $t\geq 2$, then $\mz_3^2$ is $4$PCI, contradicting Lemma~\ref{non-PCIs}. Thus, $N\cong\mz_3$. If the centralizer $C_G(N)\not=N$ then $G$ has a subgroup isomorphic to $\mz_6$, $\mz_3^2$ or $\mz_9$, which are impossible by Lemma~\ref{non-PCIs} and Lemma~\ref{pci2346}. Thus,  $C_G(N)=N$, and by the $N/C$-Theorem, $G/N$ is isomorphic to a subgroup of $\Aut(N)\cong\mz_2$. It follows that $G=\mz_3$ or $D_6$, and by Lemma~\ref{non-PCIs}, $4\leq m\leq 5$. This completes the proof of part~(1).

To prove part~(2), we assume that $G$ is $m$PDCI with $m\geq 4$. By Proposition~\ref{relation}, $G$ is $m$PCI, and by part~(1), either $G=\mz_1, \mz_2$, or $m=4,5$ and $G=\mz_3, D_6$.
We have proved $G=\mz_1$ or $\mz_2$ is $m$PDCI-group for any positive integer $m$. To finish the proof, we only need to show that $G$ is not $4$PDCI for $G=\mz_3$ or $D_6$.
First let $G=\mathbb{Z}_3=\langle x\rangle$.  Take $S_{1,2}=S_{3,2}=S_{4,1}=\{1,x,x^2\}$,
$S_{1,3}=\{x\}$,
$S_{2,4}=\{1,x\}$, $S_{3,1}=S_{4,2}=\{1\}$, and
$S_{i,j}=\emptyset$ for all other $1\leq i,j\leq 4$. With the help of {\sc Magma}~\cite{magma}, $\Gamma=\Cay(G,S_{i,j}: 1\leq i,j\leq 4)$ is a normal $4$-PCayley digraph and $\Aut(\Gamma)\cong \mathbb{Z}_3\times \mathbb{Z}_3$ has a semiregular subgroup isomorphic to $G$ with orbit set $\G=\{G_1,G_2,G_3,G_4\}$, but different from $R(G)$. By Theorem~\ref{PCI-BabaiSimilar}, $\Gamma$ is not $4$PDCI, and hence $\mz_3$ is not $4$PDCI. Now let $G=D_6=\langle x,y \ |\ x^3=y^2=1,yxy=x^2 \rangle$ and take $S_{i,j}$ as the same as above for $\mz_3$. Again by {\sc Magma}, $\Cay(D_6,S_{i,j}: 1\leq i,j\leq 4)$ is not $4$PDCI and hence $D_6$ is not $4$PDCI.
\end{proof}

\section*{Acknowledgments}
This work was partially supported by the National Natural Science Foundation of China (12331013, 12311530692, 12271024, 12161141005, 12425111,  12301461) and the 111 Project of China (B16002).
	
\section*{Conflict of interest}
The authors declare they have no financial interests.
	
\section*{Availability of data and materials}
Data sharing not applicable to this article as no datasets were generated or analysed during the current study.

{}
\end{document}